\documentclass[11pt,reqno]{amsart}
\usepackage{inputenc}
\usepackage{enumerate}
\usepackage{amssymb}
\usepackage{amsmath}
\usepackage{mathrsfs}
\usepackage{amsthm}
\usepackage{tikz-cd}
\usepackage{mathpazo}
\usepackage{extarrows}
\usepackage{color}
\usepackage{setspace}
\usepackage[colorlinks,linkcolor = blue,anchorcolor = red,urlcolor = blue,citecolor = blue]{hyperref}
\usepackage{graphicx, subfig}
\usepackage[square, comma, sort&compress, numbers]{natbib}

\theoremstyle{definition}
\newtheorem{Def}{Definition}[section]

\newtheorem{Rem}[Def]{Remark}
\theoremstyle{plain}
\newtheorem{Thm}[Def]{Theorem}
\newtheorem{Lem}[Def]{Lemma}
\newtheorem{Pro}[Def]{Proposition}
\newtheorem{Cor}[Def]{Corollary}

\newtheorem*{Que}{Question}
\usepackage[left=2.6cm, right=2.6cm, top=3cm, bottom=3cm]{geometry}
\usepackage[all]{xy}
\usepackage{bm}
\usepackage{comment}

\def\IN{\mathbb N}\def\IR{\mathbb R}\def\IC{\mathbb C}\def\IL{\mathbb L}\def\IZ{\mathbb Z}

\def\B{\mathcal B}\def\K{\mathcal K}\def\Q{\mathcal Q}\def\H{\mathcal H}\def\S{\mathcal S}\def\G{\mathcal G}

\def\supp{\textup{supp}}

\def\prop{\textup{Prop}}

\def\rank{\textup{rank}}

\def\ox{\otimes}

\def\wt{\widetilde}

\author[L.~Guo]{Liang Guo}
\address[L.~Guo]{Centre for Mathematical Sciences, Lund University, Box 118, Lund, 22100, Sweden}
\email{liang.guo@math.lth.se}

\author[J.~Qian]{Jin Qian}
\address[J.~Qian]{Research Center for Operator Algebras, School of Mathematical Sciences, East China Normal University, Shanghai, 200241, P. R. China}
\email{52265500013@stu.ecnu.edu.cn}

\author[Q.~Wang]{Qin Wang}
\address[Q.~Wang]{Research Center for Operator Algebras, 
School of Mathematical Sciences, East China Normal University, Shanghai, 200241, P. R. China.}
\email{qwang@math.ecnu.edu.cn}

\title{Ultrarank and ideal structure of uniform Roe algebras}
\date{\today}

\begin{document}

\begin{abstract}
In this paper, we study the ideal structure of the uniform Roe algebra $C^*_u(X)$ for a coarse disjoint union $X = \bigsqcup_{n\in\IN} X_n$ of finite metric spaces, with a focus on ideals lying between the compact ideal $\K$ and the ghost ideal $\G$. To analyze these ideals, we introduce a nonstandard invariant called the \emph{$\varepsilon$-ultrarank} of an operator, which is stable under compact perturbations and completely characterizes principal ideal inclusions in the block-diagonal Calkin algebra. Using the $\varepsilon$-ultrarank as an obstruction to ideal inclusion in $C^*_u(X)$, we construct an ideal between $\K$ and $\G$ that cannot be written as the intersection of any geometric ideal and the ghost ideal. On the positive side, when the averaging projection $P$ belongs to $C^*_u(X)$, we prove that the closed two-sided ideal $\langle P \rangle$ generated by $P$ equals $\G$ if and only if every ghost operator has finite ultrarank and $X$ has zero-density Property A, and we classify all closed two-sided ideals of $C^*_u(X)$ in this case.
\end{abstract}

\maketitle

\tableofcontents

\section{Introduction}

For a discrete metric space $X$ with bounded geometry, its uniform Roe algebra, denoted by $C^*_u(X)$, is the $C^*$-algebra generated by all bounded operators on $\ell^2(X)$ with finite propagation. It encodes the coarse geometric information of the underlying space $X$ and completely determines the coarse geometric structure of $X$, see \cite{BBFKVW2022}. Compared with the standard Roe algebra, the uniform Roe algebra $C^*_u(X)$ has a concrete structure that is convenient to work with, and it is closely related to many other topics, such as the Fredholm theory of band-dominated operators \cite{Roe2003, RRS2004, SpaWil2017}.

A long-standing problem in the theory of uniform Roe algebras is to classify all closed two-sided ideals of $C^*_u(X)$ in terms of the coarse geometry of $X$. This program was initiated by X.~Chen and Q.~Wang in \cite{CW2004, ChenWang2005}. Recall that an ideal $I \triangleleft C^*_u(X)$ is called \emph{geometric} if the finite propagation operators in $I$ are norm-dense in $I$. Using Skandalis-Tu-Yu coarse groupoid realization $C^*_u(X) \cong C^*_r(G(X))$  \cite{STY2002}, each invariant open subset $U \subseteq \beta X$ canonically determines a geometric ideal $I(U) = C^*_r(G(X)_U)$. Conversely, any closed two-sided ideal $J \triangleleft C^*_u(X)$ determines an invariant open subset $U(J) \subseteq \beta X$ via $J \cap \ell^\infty(X) = C_0(U(J))$. When $X$ has Yu's Property A \cite{Yu2000}, X.~Chen and Q.~Wang proved that the assignment $U \mapsto I(U)$ establishes a lattice isomorphism between the invariant open subsets of $\beta X$ and the closed two-sided ideals of $C^*_u(X)$. In particular, every ideal of $C^*_u(X)$ is geometric.

When the space $X$ lacks Property A, the ideal structure becomes substantially more intricate due to the presence of non-compact \emph{ghost operators}. By a theorem of J.~Roe and R.~Willett \cite{RoeWil2014}, $X$ fails to have Property A if and only if the ghost ideal $\G \triangleleft C^*_u(X)$, which contains all the ghost operators, strictly contains the compact ideal $\K = \K(\ell^2(X))$. Since $\K = I(X)$ and $U(\G) = X$, any intermediate ideal $\K \subsetneq J \subseteq \G$ shares the same invariant open subset $U(J) = X$.

Recently, in \cite{WZ2025}, Q.~Wang and J.~Zhang extended the ideal structure theory beyond Property A by introducing the notion of \emph{ghostly ideals}. For each invariant open subset $U \subseteq \beta X$, the ghostly ideal $\wt I(U)$ consists of all operators in $C^*_u(X)$ whose limit operators (in the sense of \v{S}pakula and Willett \cite{SpaWil2017}, or equivalently in the nonstandard limit space framework of \cite{GQW2024}) vanish on $\beta X \setminus U$. In particular, $\wt I(X) = \G$ is the classical ghost ideal. Q.~Wang and J.~Zhang proved the sandwich theorem \cite[Theorem 4.4]{WZ2025} that every ideal $J \triangleleft C^*_u(X)$ with $U(J) = U$ satisfies
$$I(U) \subseteq J \subseteq \wt I(U).$$
Observing that for any invariant open subset $V \supseteq U$, the intersection $I(V) \cap \wt I(U)$ lies in the fiber $\mathfrak{I}_U = \{J \triangleleft C^*_u(X) \mid U(J) = U\}$, Q.~Wang and J.~Zhang posed the following open question (see \cite[Question 8.1]{WZ2025}):

\begin{Que}
Let $X$ be a discrete metric space with bounded geometry. Is it true that for every closed two-sided ideal $J \triangleleft C^*_u(X)$, there exists an invariant open subset $V \supseteq U(J)$ in $\beta X$ such that
\begin{equation}\label{eq: WZ question}
J = I(V) \cap \wt I(U(J))\,?
\end{equation}
\end{Que}

An affirmative answer to this question would yield a complete geometric and ghostly classification of the ideal lattice of $C^*_u(X)$. The main purpose of this paper is to investigate this question for sparse spaces $X = \bigsqcup_{n\in\IN} X_n$ by introducing a new nonstandard invariant, which we call the \emph{$\varepsilon$-ultrarank}.

In \cite{GQW2024}, the authors developed an ultraproduct approach to limit space theory by fixing a non-principal ultrafilter $\omega \in \partial_\beta \IN$ and equipping the ultrapower $X^\omega$ with a hyperreal-valued metric $d^\omega: X^\omega \times X^\omega \to \IR^\omega$. In the present paper, for a coarse disjoint union of finite metric spaces $X = \bigsqcup_{n\in\IN} X_n$ with $d(X_n, X \setminus X_n) \to \infty$, any operator $T \in C^*_u(X)$ is block-diagonal modulo $\K$, i.e., $T - \bigoplus_{n\in\IN} T_n \in \K$, where $T_n = \chi_{X_n} T \chi_{X_n}$. For each $\varepsilon > 0$, we define the $\varepsilon$-rank $r_\varepsilon(T_n)$ to be the number of singular values of $T_n$ strictly greater than $\varepsilon$, and define the \emph{$\varepsilon$-ultrarank} of $T$ along $\omega$ to be the hypernatural number
$$\rho_{\varepsilon,\omega}(T) = [r_\varepsilon(T_n)]_\omega \in \IN^\omega \subseteq \IR^\omega.$$
The ultrarank profile $\{\rho_{\varepsilon,\omega}(T)\}_{\varepsilon > 0}$ is robust under norm approximations and compact perturbations. We write $T \preceq_\omega S$ if for every $\varepsilon > 0$, there exist $\delta > 0$ and an integer $N \ge 1$ such that $\rho_{\varepsilon,\omega}(T) \le N\, \rho_{\delta,\omega}(S)$, and write $T \preceq S$ if $T \preceq_\omega S$ holds for all $\omega \in \partial_\beta \IN$.

For a $C^*$-algebra $A$ and an element $S \in A$, we denote by $\langle S \rangle_A$ the closed two-sided ideal of $A$ generated by $S$, omitting the subscript when $A = C^*_u(X)$. We first show that in the block-diagonal Calkin algebra $\Q = \prod_{n\in\IN} \B(\ell^2(X_n)) / \bigoplus_{n\in\IN} \B(\ell^2(X_n))$, where no finite propagation constraint is imposed, the ultrarank order $T \preceq S$ is both necessary and sufficient for $[T] \in \langle [S] \rangle_{\Q}$ (Theorem \ref{thm: ambient ideal ultrarank}). Passing to $C^*_u(X)$, the order $T \preceq S$ provides an obstruction to ideal inclusion $T \in \langle S \rangle$ (Corollary \ref{cor: ultrarank obstruction}): if $T \not\preceq_\omega S$ for some $\omega \in \partial_\beta \IN$, then $T$ cannot belong to $\langle S \rangle$. By considering the coarse disjoint union $X = \bigsqcup_{n\in\IN} (E_n \times E_n)$ of Cartesian products of expander graphs, we use this obstruction to obtain a negative answer to \cite[Question 8.1]{WZ2025} in general:

\begin{Thm}[Theorem \ref{thm: counterexample WZ}]
Let $(E_n)_{n\in\IN}$ be a sequence of finite connected $d$-regular expander graphs, and let $X = \bigsqcup_{n\in\IN} (E_n \times E_n)$ be their coarse disjoint union equipped with the $\ell^1$-product metric on each block. Let $J_P = \langle P \rangle$ be the ideal of $C^*_u(X)$ generated by the averaging projection $P = \bigoplus_{n\in\IN} (p_n \ox p_n)$. Then $\K \subsetneq J_P \subsetneq \G$ and $U(J_P) = X$, but there does not exist any geometric ideal $I \triangleleft C^*_u(X)$ such that
$$J_P = I \cap \G.$$
\end{Thm}

Given that $T \preceq S$ is necessary for $T \in \langle S \rangle$ in $C^*_u(X)$ and sufficient in the block-diagonal Calkin algebra $\Q$, it is natural to ask whether the converse holds in $C^*_u(X)$: if $T \preceq S$ (or even $T \prec_\omega S$ for all $\omega \in \partial_\beta \IN$), does $\langle S \rangle$ necessarily contain $T$? We show that this fails in $C^*_u(X)$. Intuitively, the geometric regions supporting the singular vectors of $T_n$ and $S_n$ may drift infinitely far apart as $n \to \infty$, preventing any finite propagation operator from intertwining them.

\begin{Thm}[Theorem \ref{thm: low ultrarank no absorption}]
In the product-expander space $X = \bigsqcup_{n\in\IN} (E_n \times E_n)$, let $P_1 = \bigoplus_{n\in\IN} (p_n \ox 1)$ be the first-coordinate averaging ghost projection and $P_\Delta$ be the rank-one averaging ghost projection supported on the diagonal $\Delta_n = \{(a,a) \mid a \in E_n\} \subseteq E_n \times E_n$. Then for every $\varepsilon \in (0,1)$ and every non-principal ultrafilter $\omega \in \partial_\beta \IN$,
$$\rho_{\varepsilon,\omega}(P_\Delta) = 1 \quad \text{and} \quad \rho_{\varepsilon,\omega}(P_1) = [\# E_n]_\omega,$$
so that $P_\Delta \prec_\omega P_1$. Nevertheless, $P_1 C^*_u(X) P_\Delta \subseteq \K$, and consequently
$$\langle P_\Delta \rangle \cap \langle P_1 \rangle = \K \quad \text{and} \quad P_\Delta \notin \langle P_1 \rangle.$$
\end{Thm}

In the final part of the paper, we turn to positive results for \cite[Question 8.1]{WZ2025}. For a coarse disjoint union $X = \bigsqcup_{n\in\IN} X_n$ with $\# X_n \to \infty$, the averaging projection $P = \bigoplus_{n\in\IN} p_n$ belongs to $C^*_u(X)$ if and only if $(X_n)_{n\in\IN}$ is an asymptotic expander \cite{KLVZ2021}. When $P \in C^*_u(X)$, we first show that an operator $T \in C^*_u(X)$ belongs to $\langle P \rangle$ if and only if $T \preceq P$ and $T$ is \emph{uniformly continuous in density} (Proposition \ref{thm: characterization of JP}). Using a theorem of J.~Roe and R.~Willett \cite{RoeWil2014}, we relate uniform continuity in density for ghost operators to Yu's Property A on zero-density subspaces, which yields a complete characterization of $\G = \langle P \rangle$ and an affirmative answer to \cite[Question 8.1]{WZ2025} when $\G = \langle P \rangle$:

\begin{Thm}[Theorem \ref{thm: G=JP characterization} and Corollary \ref{cor: general ideal classification}]
Let $X = \bigsqcup_{n\in\IN} X_n$ be the coarse disjoint union of a sequence of asymptotic expanders. Then:
\begin{itemize}
\item[(1)] $P$ generates the ghost ideal $\G = \langle P \rangle$ if and only if both of the following conditions hold:
\begin{itemize}
\item[(i)] \emph{(Finite ultrarank)} For every $T \in \G$, one has $T \preceq P$, equivalently, $\sup_{n\in\IN} r_\varepsilon(T_n) < \infty$;
\item[(ii)] \emph{(Zero-density Property A)} For every sequence of subsets $Y_n \subseteq X_n$ satisfying $\lim_{n\to\infty} \frac{\# Y_n}{\# X_n} = 0$, the subspace $Y = \bigsqcup_{n\in\IN} Y_n$ has Yu's Property A.
\end{itemize}
\item[(2)] If $\G = \langle P \rangle$, then for every nonzero closed two-sided ideal $J \triangleleft C^*_u(X)$, there exists an invariant open subset $V \supseteq U(J)$ in $\beta X$ such that
$$J = I(V) \cap \wt I(U(J)).$$
\end{itemize}
\end{Thm}

The paper is organized as follows. In Section \ref{sec: ultrarank}, we review the background on uniform Roe algebras, geometric and ghostly ideals, and then introduce the $\varepsilon$-ultrarank of operators and the order $\preceq$. In Section \ref{sec: counterexample}, we prove the ultrarank characterization in $\Q$ and the obstruction for ideal inclusion in $C^*_u(X)$, and construct the product-expander counterexample to \cite[Question 8.1]{WZ2025}. In Section \ref{sec: converse failure}, we prove Theorem \ref{thm: low ultrarank no absorption}, showing that lower-order $\varepsilon$-ultrarank does not imply ideal inclusion in $C^*_u(X)$. In Section \ref{sec: positive results}, we establish the operator and ultrarank characterizations of $\langle P \rangle$ and $\G = \langle P \rangle$, and prove the positive classification results for \cite[Question 8.1]{WZ2025} in the case where $\G = \langle P \rangle$.

\subsection*{Acknowledgment}
This work is partially supported by the National Natural Science Foundation of China (NSFC), Grant No. 12571135.
L.~Guo is supported by the Swedish Foundation for International Cooperation in Research and Higher Education (STINT) (Grant No. PS2025-9545).

\subsection*{AI declaration}
During the preparation of this paper, we initially used GPT 6.1 Sol to find the example in Theorem~\ref{thm: counterexample WZ}. From this example, the authors formulated the notion of ultrarank and independently developed its theory, including the characterization of principal ideals without propagation constraints (Theorem~\ref{thm: ambient ideal ultrarank}) and its failure in $C^*_u(X)$ (Theorem~\ref{thm: low ultrarank no absorption}). In Section~\ref{sec: positive results}, AI tools were used to assist in completing the proof of Lemma~\ref{lem: ghostly decomposition}. AI tools were also used for language polishing, translation, and improving the presentation of the proofs. All remaining results and proofs are the authors' own work, and the authors have carefully verified all arguments and bear full responsibility for the manuscript.

\section{Preliminaries and the ultrarank of operators}\label{sec: ultrarank}

In this section, we shall first recall the necessary background on uniform Roe algebras, geometric and ghostly ideals \cite{CW2004, WZ2025}, and the nonstandard ultraproduct framework from \cite{GQW2024}. Based on these, we introduce a new invariant called the \emph{$\varepsilon$-ultrarank} for operators on sparse metric spaces.

\subsection{Ideals of uniform Roe algebras}

Throughout this paper, let $(X, d)$ be a countable, strongly discrete metric space with bounded geometry. Recall from \cite{SpaWil2017} that $X$ is \emph{strongly discrete} if $d(x,y) \in \IZ$ for all $x, y \in X$, and $X$ has \emph{bounded geometry} if for every $R > 0$,
$$N_R := \sup_{x \in X} \# B(x, R) < \infty.$$
For any bounded operator $T \in \B(\ell^2(X))$ and $x, y \in X$, the $(x,y)$-th matrix entry of $T$ is written as $T(x,y) = \langle \delta_x, T \delta_y \rangle$. The \emph{propagation} of $T$ is defined by
$$\prop(T) = \sup\{ d(x,y) \mid x, y \in X,\ T(x,y) \neq 0 \}.$$
The algebraic uniform Roe algebra $\IC_u[X]$ consists of all operators $T \in \B(\ell^2(X))$ with finite propagation, and the \emph{uniform Roe algebra} $C^*_u(X)$ is the operator norm closure of $\IC_u[X]$ in $\B(\ell^2(X))$. Denote by $\K = \K(\ell^2(X))$ the ideal of compact operators. Throughout this paper, for a $C^*$-algebra $A$ and a subset (or element) $\S \subseteq A$, we denote by $\langle \S \rangle_A$ the closed two-sided ideal of $A$ generated by $\S$, and simply write $\langle \S \rangle$ when $A = C^*_u(X)$.

\begin{Def}
An operator $T \in C^*_u(X)$ is called a \emph{ghost operator} if for every $\varepsilon > 0$, there exists a finite subset $F \subseteq X$ such that $|T(x,y)| < \varepsilon$ whenever $(x,y) \notin F \times F$. The set of all ghost operators forms a closed two-sided ideal in $C^*_u(X)$, called the \emph{ghost ideal} and denoted by $\G$.
\end{Def}

To describe the general ideal structure of $C^*_u(X)$, we recall the notation of Chen-Wang \cite{CW2004} and Wang-Zhang \cite{WZ2025}. Identify $\ell^\infty(X)$ with $C(\beta X)$, where $\beta X$ is the Stone-\v{C}ech compactification of $X$, and write $\partial_\beta X = \beta X \setminus X$ for the Stone-\v{C}ech boundary. For a subset $B \subseteq X$, we denote by $\chi_B \in \ell^\infty(X) \subseteq C^*_u(X)$ the diagonal projection onto $\ell^2(B)$. For a closed two-sided ideal $J \triangleleft C^*_u(X)$, its intersection $J \cap \ell^\infty(X)$ is a closed ideal of $C(\beta X)$, hence is of the form $C_0(U(J))$ for an open subset $U(J) \subseteq \beta X$. The open set $U(J)$ is invariant under the coarse groupoid $G(X)$ (equivalently, closed under compatible partial translations, see \cite[Section 3]{WZ2025}).

Conversely, let $U \subseteq \beta X$ be an invariant open set. For any operator $T \in \B(\ell^2(X))$ and $\varepsilon > 0$, define the $\varepsilon$-support of $T$ by
$$\supp_\varepsilon(T) = \{ (x,y) \in X \times X \mid |T(x,y)| \ge \varepsilon \},$$
and let $r: G(X) \to \beta X$ denote the range map of the coarse groupoid (whose restriction to $X \times X$ is the first-coordinate projection $(x,y) \mapsto x$). There are two canonical ideals of $C^*_u(X)$ associated with $U$:
\begin{itemize}
\item[(1)] The \emph{geometric ideal} of $C^*_u(X)$ with respect to $U$, denoted by $I(U)$, is defined to be the norm closure of
$$\IL[U] = \left\{ T \in \IC_u[X]\ \big|\ \overline{r(\supp_\varepsilon(T))}^{\beta X} \subseteq U \text{ for every } \varepsilon > 0 \right\}.$$
Equivalently, $I(U)$ is the closed ideal of $C^*_u(X)$ generated by $C_0(U) \subseteq \ell^\infty(X)$, that is,
$$I(U) = \big\langle \{ \chi_B \mid B \subseteq X,\ \overline{B}^{\beta X} \subseteq U \} \big\rangle.$$
An ideal $J \triangleleft C^*_u(X)$ is geometric (i.e., $J = \overline{J \cap \IC_u[X]}^{\|\cdot\|}$) if and only if $J = I(U(J))$. In particular, $I(X) = \K$ and $I(\beta X) = C^*_u(X)$.
\item[(2)] The \emph{ghostly ideal} of $C^*_u(X)$ with respect to $U$, denoted by $\wt I(U)$, is defined to be
$$\wt I(U) = \left\{ T \in C^*_u(X)\ \big|\ \overline{r(\supp_\varepsilon(T))}^{\beta X} \subseteq U \text{ for every } \varepsilon > 0 \right\}.$$
Equivalently, $\wt I(U)$ consists of all operators $T \in C^*_u(X)$ whose limit operators vanish outside $U$. In particular, $\wt I(X) = \G$ and $\wt I(\beta X) = C^*_u(X)$.
\end{itemize}

For every invariant open subset $U \subseteq \beta X$, one has $U(I(U)) = U(\wt I(U)) = U$, and every ideal $J \triangleleft C^*_u(X)$ with $U(J) = U$ satisfies the sandwich relation \cite[Theorem 4.4]{WZ2025}:
\begin{equation}\label{eq: sandwich}
I(U) \subseteq J \subseteq \wt I(U).
\end{equation}

\subsection{Ultrarank}

From now on, we focus on a \emph{coarse disjoint union} of a sequence of finite metric spaces:
\begin{equation}\label{eq: sparse X}
X = \bigsqcup_{n\in\IN} X_n, \quad \text{with} \quad d(X_n, X \setminus X_n) \to \infty \text{ as } n \to \infty,
\end{equation}
and assume that $X$ is strongly discrete with bounded geometry. Write $\H_n = \ell^2(X_n)$. For any $T \in C^*_u(X)$, denote its $n$-th block compression and block-diagonal expectation by
$$T_n = \chi_{X_n} T \chi_{X_n} \in \B(\H_n), \quad E(T) = \bigoplus_{n\in\IN} T_n \in \B(\ell^2(X)).$$
It is a well-known fact (see, for example, \cite[Lemma 12.5.4]{WilYu2020}) that on such a coarse disjoint union, $T-E(T) \in \K$ and $\|[T]\|_{C^*_u(X)/\K} = \limsup\limits_{n\to\infty} \|T_n\|$.

We now recall the ultrapower construction from \cite[Section 2.1]{GQW2024}. Fix a non-principal ultrafilter $\omega \in \partial_\beta \IN$. For two sequences $(r_n)_{n\in\IN}, (s_n)_{n\in\IN} \in \prod_{n\in\IN} \IR$, define the equivalence relation $(r_n) \sim_{s,\omega} (s_n)$ by
$$\{ n \in \IN \mid r_n = s_n \} \in \omega.$$
Let $\IR^\omega = \prod_{n\in\IN} \IR / \sim_{s,\omega}$ be the hyperreal ordered field, and let $\IN^\omega \subseteq \IR^\omega$ be the set of hypernatural numbers (internal hyperintegers), consisting of $\omega$-equivalence classes $[k_n]_\omega$ of sequences $(k_n)_{n\in\IN}$ in $\IN \cup \{0\}$. As in \cite[Section 2.1]{GQW2024}, an element $[r_n]_\omega \in \IR^\omega$ is called \emph{finite} if there exists a standard real number $M > 0$ such that $[|r_n|]_\omega \le M$, \emph{infinite} if $[|r_n|]_\omega > M$ for all $M > 0$, and \emph{infinitesimal} if $[|r_n|]_\omega < \varepsilon$ for all $\varepsilon > 0$. Notice that a hypernatural number $[k_n]_\omega \in \IN^\omega$ is finite if and only if it is equal to a standard integer $N \in \IN \cup \{0\}$ (since $\omega$ is an ultrafilter and $\{0,1,\dots,M\}$ is a finite partition).

\begin{Def}[$\varepsilon$-rank and $\varepsilon$-ultrarank]\label{def: ultrarank}
Let $T \in \B(\H_n)$ be a matrix on the finite-dimensional Hilbert space $\H_n$, and let $s_1(T) \ge s_2(T) \ge \cdots \ge s_{\# X_n}(T) \ge 0$ be its singular values (i.e., the eigenvalues of $|T| = (T^*T)^{1/2}$) repeated according to multiplicity. By the singular value decomposition (SVD, see e.g.~\cite[Theorem 2.6.3]{HJ2013}), there exist orthonormal bases $\{v_j\}_{j=1}^{\# X_n}$ and $\{u_j\}_{j=1}^{\# X_n}$ of $\H_n$ such that
$$T = \sum_{j=1}^{\# X_n} s_j(T)\, (u_j \ox v_j),$$
where $u \ox v \in \B(\H_n)$ denotes the rank-one operator defined by $(u \ox v)(\xi) = \langle v, \xi \rangle u$ for $\xi \in \H_n$.
\begin{itemize}
\item[(1)] For any $\varepsilon > 0$, the \emph{$\varepsilon$-rank} of $T$ is defined by
\begin{equation}\label{eq: r_epsilon}
r_\varepsilon(T) := \#\{ j \in \{1,\dots,\# X_n\} \mid s_j(T) > \varepsilon \}.
\end{equation}
\item[(2)] Let $\omega \in \partial_\beta \IN$ be a non-principal ultrafilter. For any $T \in C^*_u(X)$ (or any bounded sequence $(T_n) \in \prod_{n\in\IN} \B(\H_n)$) and any $\varepsilon > 0$, the \emph{$\varepsilon$-ultrarank} of $T$ with respect to $\omega$ is defined to be the hypernatural number
\begin{equation}\label{eq: ultrarank}
\rho_{\varepsilon,\omega}(T) := [r_\varepsilon(T_n)]_\omega \in \IN^\omega,
\end{equation}
where $T_n = \chi_{X_n} T \chi_{X_n}$. The family $\{\rho_{\varepsilon,\omega}(T)\}_{\varepsilon > 0}$ is called the \emph{ultrarank profile} of $T$.
\end{itemize}
\end{Def}

The basic properties of the $\varepsilon$-rank and the $\varepsilon$-ultrarank are summarized in the following.

\begin{Lem}\label{lem: basic rank inequalities}
For any finite-dimensional matrices $a, b, c \in \B(\H_n)$ and $\varepsilon > 0$, one has that
\begin{itemize}
\item[(1)] $r_\varepsilon(a) = \min\{ \rank(a') \mid a' \in \B(\H_n),\ \|a - a'\| \le \varepsilon \}$;
\item[(2)] $r_\varepsilon(a + b) \le r_{\varepsilon/2}(a) + r_{\varepsilon/2}(b)$;
\item[(3)] $r_\varepsilon(cab) \le r_{\varepsilon / (\|c\|\|b\|)}(a)$ for all $c,b \ne 0$;
\item[(4)] if $\|a - b\| < \varepsilon - \eta$ for some $0 < \eta < \varepsilon$, then $r_\varepsilon(a) \le r_\eta(b)$.
\end{itemize}
\end{Lem}

\begin{proof}
To prove (1), write the SVD (see \cite[Theorem 2.6.3]{HJ2013}) of $a$ as $a = \sum_{j=1}^{\# X_n} s_j(a)\, (u_j \ox v_j)$ with two orthonormal bases $\{u_j\}$ and $\{v_j\}$, and set $k = r_\varepsilon(a)$. Taking the truncated operator $a' = \sum_{j=1}^k s_j(a)\, (u_j \ox v_j)$, we have $\rank(a') = k$ and $\|a - a'\| = s_{k+1}(a) \le \varepsilon$, which shows that the minimum is at most $k$. Conversely, suppose that $a' \in \B(\H_n)$ satisfies $\|a - a'\| \le \varepsilon$ and $\rank(a') < k$. Let $W = \textup{span}\{v_1, \dots, v_k\} \subseteq \H_n$. Since $\dim(W) = k$ and $\dim(\ker(a')) = \# X_n - \rank(a') > \# X_n - k$, we have $\dim(W) + \dim(\ker(a')) > \# X_n$, so $W \cap \ker(a') \neq \{0\}$. Choose a unit vector $v = \sum_{j=1}^k c_j v_j \in W \cap \ker(a')$ with $\|v\| = 1$. Since $a'v = 0$ and $s_j(a) > \varepsilon$ for all $1 \le j \le k$, we obtain
$$\|a - a'\| \ge \|(a - a')v\| = \|av\| = \Big\| \sum_{j=1}^k c_j s_j(a) u_j \Big\| = \Big( \sum_{j=1}^k |c_j|^2 s_j(a)^2 \Big)^{1/2} > \varepsilon,$$
a contradiction. Thus $\rank(a') \ge k = r_\varepsilon(a)$.

Items (2), (3), and (4) follow immediately from (1) by taking optimal low-rank approximations for $a$ and $b$.
\end{proof}

\begin{Def}[Order of ultrarank profiles]\label{def: ultrarank order}
Let $T = (T_n)$ and $S = (S_n)$ be bounded sequences in $\prod_{n\in\IN} \B(\H_n)$.
\begin{itemize}
\item[(1)] Given a non-principal ultrafilter $\omega \in \partial_\beta \IN$, we say that the ultrarank of $T$ is of \emph{at most the same order as} that of $S$ along $\omega$, denoted by $T \preceq_\omega S$, if for every $\varepsilon > 0$, there exist $\delta > 0$ and an integer $N \ge 1$ such that
$$\rho_{\varepsilon,\omega}(T) \le N\, \rho_{\delta,\omega}(S).$$
We write $T \simeq_\omega S$ if both $T \preceq_\omega S$ and $S \preceq_\omega T$ hold, and write $T \prec_\omega S$ if $T \preceq_\omega S$ but $S \not\preceq_\omega T$.
\item[(2)] We write $T \preceq S$ if $T \preceq_\omega S$ holds for all $\omega \in \partial_\beta \IN$, and write $T \simeq S$ if $T \preceq S$ and $S \preceq T$.
\end{itemize}
\end{Def}

It is routine to check that $\preceq_\omega$ and $\preceq$ are reflexive and transitive, so they define partial orders on the $\simeq_\omega$- and $\simeq$-equivalence classes of operators, respectively. The following proposition shows that the ultrarank profile is stable under compact perturbations, and hence its $\simeq$-equivalence class is invariant in the Calkin algebra $C^*_u(X)/\K$.

\begin{Pro}\label{pro: ultrarank stability}
Let $T, S \in C^*_u(X)$ be such that $T - S \in \K$. Then for any $0 < \eta < \varepsilon$ and any non-principal ultrafilter $\omega \in \partial_\beta \IN$, one has that
$$\rho_{\varepsilon,\omega}(T) \le \rho_{\eta,\omega}(S) \quad \text{and} \quad \rho_{\varepsilon,\omega}(S) \le \rho_{\eta,\omega}(T).$$
Consequently, $T \simeq S$. In particular, since $T - E(T) \in \K$, one has $T \simeq E(T)$.
\end{Pro}

\begin{proof}
Since $T - S \in \K$, we have $\lim_{n\to\infty} \|T_n - S_n\| = 0$. Thus $\|T_n - S_n\| < \varepsilon - \eta$ for all sufficiently large $n \in \IN$. Applying Lemma \ref{lem: basic rank inequalities} (4) gives $r_\varepsilon(T_n) \le r_\eta(S_n)$ and $r_\varepsilon(S_n) \le r_\eta(T_n)$ on a cofinite subset of $\IN$, which belongs to $\omega$. Taking $\eta = \varepsilon/2$ and $N = 1$ in Definition \ref{def: ultrarank order} yields $T \simeq S$.
\end{proof}

\begin{Rem}\label{rem: why not exact rank}
In view of Lemma \ref{lem: basic rank inequalities} (1), $r_\varepsilon(T_n)$ is the minimal rank of an operator required to approximate $T_n$ within $\varepsilon$ in norm. Thus, the order $T \preceq_\omega S$ admits an interpretation in terms of low-rank approximations: it means that for any given error $\varepsilon > 0$, the minimal rank growth needed to $\varepsilon$-approximate $T_n$ along $\omega$ can be controlled (up to a finite multiple $N$) by the approximation rank growth of $S_n$ at some precision $\delta > 0$. In particular, $T \simeq_\omega S$ means that $T_n$ and $S_n$ share the same asymptotic rate of rank growth under norm approximations.

This approximation perspective also explains why it is essential to work with the full ultrarank profile $\{\rho_{\varepsilon,\omega}(T)\}_{\varepsilon > 0}$ up to $\simeq$, rather than the exact rank $[\rank(T_n)]_\omega$ or a single $\varepsilon > 0$. For instance, the compact operator $C = \bigoplus_{n\in\IN} \frac{1}{n+1} I_{\H_n} \in \K$ has maximal exact rank $[\rank(C_n)]_\omega = [\# X_n]_\omega$, yet for any $\varepsilon > 0$ it can be $\varepsilon$-approximated by $0$ for large $n$, so that $C \simeq 0$. Similarly, at a single fixed $\varepsilon > 0$, a compact perturbation can shift singular values across $\varepsilon$ (compare $S = \varepsilon I$ with $T = S + C$, where $\rho_{\varepsilon,\omega}(S) = 0$ while $\rho_{\varepsilon,\omega}(T) = [\# X_n]_\omega$), which is absorbed once we allow the precision $\delta$ to vary in Definition \ref{def: ultrarank order}.
\end{Rem}

\section{Counterexample via ultrarank}\label{sec: counterexample}

In this section, we prove that the $\varepsilon$-ultrarank provides a necessary condition for ideal inclusion in $C^*_u(X)$. Using this obstruction, we show that the Cartesian product of expander graphs gives a counterexample to \cite[Question 8.1]{WZ2025}.

\subsection{The ultrarank obstruction to ideal inclusion}

Recall that the block compression map
$$E: C^*_u(X)/\K \to \Q = \prod_{n\in\IN} \B(\H_n) / \bigoplus_{n\in\IN} \B(\H_n), \quad E(T + \K) = [(T_n)],$$
is a faithful $*$-homomorphism. Before turning to the ideals of $C^*_u(X)$, we first investigate principal ideal inclusions in the block-diagonal Calkin algebra $\Q$, where no propagation constraint is imposed on the operators.

This consideration is motivated by the classical ideal structure of $\B(\H)$. It is a well-known fact that for a separable Hilbert space $\H$, the compact ideal $\K(\H)$ is the unique non-trivial closed two-sided ideal of $\B(\H)$. When $\H$ is non-separable, however, operators with countable range generate proper intermediate ideals, and the ideal lattice of $\B(\H)$ is completely classified by the infinite cardinalities of the $\varepsilon$-ranks of operators. In our setting of a coarse disjoint union $X = \bigsqcup_{n\in\IN} X_n$, although each block $\H_n$ is finite-dimensional, the growth rate of the $\varepsilon$-rank $r_\varepsilon(T_n)$ as $n \to \infty$ (measured by the $\varepsilon$-ultrarank $\rho_{\varepsilon,\omega}(T) \in \IN^\omega$) plays a role analogous to the rank cardinality on a non-separable Hilbert space. 

As the following theorem shows, if we drop the finite propagation condition and work inside $\Q$, the growth rate of the $\varepsilon$-ultrarank completely characterizes when one principal ideal is contained in another.

\begin{Thm}\label{thm: ambient ideal ultrarank}
Let $S = (S_n)$ and $T = (T_n)$ be bounded sequences in $\prod_{n\in\IN} \B(\H_n)$. The following are equivalent:
\begin{itemize}
\item[(1)] $[T]$ belongs to the closed two-sided ideal $\langle [S] \rangle_{\Q}$ generated by $[S]$ in $\Q$;
\item[(2)] for every $\varepsilon > 0$, there exist $\delta > 0$, a standard integer $N \ge 1$, and $n_0 \in \IN$ such that
\begin{equation}\label{eq: tail rank domination}
r_\varepsilon(T_n) \le N\, r_\delta(S_n) \quad \text{for all } n \ge n_0;
\end{equation}
\item[(3)] $T \preceq S$, that is, $T \preceq_\omega S$ for every non-principal ultrafilter $\omega \in \partial_\beta \IN$.
\end{itemize}
In particular, $\langle [T] \rangle_{\Q} = \langle [S] \rangle_{\Q}$ if and only if $T \simeq S$.
\end{Thm}

\begin{proof}
$(1) \Rightarrow (2)$: Suppose $[T] \in \langle [S] \rangle_{\Q}$ and fix $\varepsilon > 0$. There exist finitely many bounded sequences $A^{(i)} = (A^{(i)}_n)$ and $B^{(i)} = (B^{(i)}_n)$ ($1 \le i \le N$) in $\prod_{n\in\IN} \B(\H_n)$ such that
$$\Big\| [T] - \sum_{i=1}^N [A^{(i)}][S][B^{(i)}] \Big\|_\Q < \frac{\varepsilon}{4}.$$
Choose $\delta > 0$ small enough so that $\delta \sum_{i=1}^N \|A^{(i)}\| \|B^{(i)}\| < \varepsilon/4$. For each $n \in \IN$, let $S_n^\delta$ be the SVD truncation of $S_n$ retaining all singular values strictly greater than $\delta$. Then $\rank(S_n^\delta) = r_\delta(S_n)$ and $\|S_n - S_n^\delta\| \le \delta$. Set $S'_n = \sum_{i=1}^N A^{(i)}_n S_n^\delta B^{(i)}_n$. Then
$$\limsup_{n\to\infty} \|T_n - S'_n\|\leq \limsup_{n\to\infty} \Big\| T_n - \sum_{i=1}^N A^{(i)}_nS_nB^{(i)}_n \Big\| + \Big\| \sum_{i=1}^N A^{(i)}_n(S_n-S'_n)B^{(i)}_n \Big\| < \frac{\varepsilon}{4}+\frac{\varepsilon}{4} = \frac{\varepsilon}{2}.$$
Thus, there exists $n_0 \in \IN$ such that $\|T_n - S'_n\| < \varepsilon$ for all $n \ge n_0$. Since $\rank(S'_n) \le N\, r_\delta(S_n)$, Lemma \ref{lem: basic rank inequalities} (1) yields $r_\varepsilon(T_n) \le N\, r_\delta(S_n)$ for all $n \ge n_0$.

$(2) \Rightarrow (1)$: Fix $\varepsilon > 0$. By (2), there exist $\delta > 0$, an integer $N \ge 1$, and $n_0 \in \IN$ such that $r_\varepsilon(T_n) \le N\, r_\delta(S_n)$ for all $n \ge n_0$. For each $n \in \IN$, write the singular value decompositions of $T_n$ and $S_n$ as
$$T_n = \sum_{m=1}^{\# X_n} \tau_n^{(m)} \big(u_n^{(m)} \ox v_n^{(m)}\big) \quad \text{and} \quad S_n = \sum_{m=1}^{\# X_n} \sigma_n^{(m)} \big(\wt u_n^{(m)} \ox \wt v_n^{(m)}\big),$$
where $\tau_n^{(1)} \ge \cdots \ge \tau_n^{(\# X_n)} \ge 0$ and $\sigma_n^{(1)} \ge \cdots \ge \sigma_n^{(\# X_n)} \ge 0$ are the singular values of $T_n$ and $S_n$, respectively, and $\{u_n^{(m)}\}_{m=1}^{\# X_n}$, $\{v_n^{(m)}\}_{m=1}^{\# X_n}$, $\{\wt u_n^{(m)}\}_{m=1}^{\# X_n}$, $\{\wt v_n^{(m)}\}_{m=1}^{\# X_n}$ are orthonormal bases of $\H_n$. Set $p_n = r_\varepsilon(T_n)$ and $q_n = r_\delta(S_n)$, so that $\tau_n^{(m)} > \varepsilon$ for $1 \le m \le p_n$ and $\sigma_n^{(m)} > \delta$ for $1 \le m \le q_n$. Let
$$T_n^\varepsilon = \sum_{m=1}^{p_n} \tau_n^{(m)} \big(u_n^{(m)} \ox v_n^{(m)}\big)$$
be the SVD truncation of $T_n$ at $\varepsilon$, which satisfies $\|T_n - T_n^\varepsilon\| \le \varepsilon$. For convenience, we set $\tau_n^{(k)} = 0$ and $u_n^{(k)} = v_n^{(k)} = 0$ whenever $k > p_n$.

For each $n \ge n_0$, since $p_n \le N q_n$, we can partition the $p_n$ singular value components of $T_n^\varepsilon$ into $N$ groups, each of size at most $q_n$. For the $i$-th group ($1 \le i \le N$), let $B_n^{(i)} \in \B(\H_n)$ be the partial isometry mapping the right singular vectors $\{v_n^{((i-1)q_n + m)}\}_{m=1}^{q_n}$ of $T_n^\varepsilon$ in this group onto the first $q_n$ right singular vectors $\{\wt v_n^{(m)}\}_{m=1}^{q_n}$ of $S_n$, and let $A_n^{(i)} \in \B(\H_n)$ be the weighted partial isometry mapping the corresponding left singular vectors $\{\wt u_n^{(m)}\}_{m=1}^{q_n}$ of $S_n$ to $\{u_n^{((i-1)q_n + m)}\}_{m=1}^{q_n}$ scaled by the ratios $\tau_n^{((i-1)q_n + m)} / \sigma_n^{(m)}$. Explicitly, for each $1 \le i \le N$, $A_n^{(i)}$ and $B_n^{(i)}$ are given by
$$A_n^{(i)} = \sum_{m=1}^{q_n} \frac{\tau_n^{((i-1)q_n + m)}}{\sigma_n^{(m)}} \big( u_n^{((i-1)q_n + m)} \ox \wt u_n^{(m)} \big), \qquad B_n^{(i)} = \sum_{m=1}^{q_n} \big( \wt v_n^{(m)} \ox v_n^{((i-1)q_n + m)} \big)$$
for $n \ge n_0$, and we set $A_n^{(i)} = B_n^{(i)} = 0$ for $n < n_0$. Since $0 \le \tau_n^{((i-1)q_n + m)} / \sigma_n^{(m)} \le \|T\|/\delta$, we have $\|A_n^{(i)}\| \le \|T\|/\delta$ and $\|B_n^{(i)}\| \le 1$ for all $n \in \IN$. Moreover, since $S_n \wt v_n^{(m)} = \sigma_n^{(m)} \wt u_n^{(m)}$, it follows that
$$\sum_{i=1}^N A_n^{(i)} S_n B_n^{(i)} = \sum_{i=1}^N \sum_{m=1}^{q_n} \tau_n^{((i-1)q_n + m)} \big( u_n^{((i-1)q_n + m)} \ox v_n^{((i-1)q_n + m)} \big) = T_n^\varepsilon \quad (n \ge n_0).$$
Therefore, $[T^\varepsilon] = \sum_{i=1}^N [A^{(i)}][S][B^{(i)}] \in \langle [S] \rangle_{\Q}$. Since $\|[T] - [T^\varepsilon]\|_\Q \le \varepsilon$ and $\varepsilon > 0$ is arbitrary, we conclude that $[T] \in \langle [S] \rangle_{\Q}$.

$(2) \Leftrightarrow (3)$: If (2) holds, then for every $\omega \in \partial_\beta \IN$, the cofinite set $\{n \in \IN \mid n \ge n_0\}$ belongs to $\omega$, so $\rho_{\varepsilon,\omega}(T) \le N\, \rho_{\delta,\omega}(S)$ in $\IN^\omega$, which means $T \preceq S$. Conversely, suppose that (2) fails for some $\varepsilon > 0$. Then for every $k \in \IN$, the set $F_k = \{ n \in \IN \mid r_\varepsilon(T_n) > k\, r_{1/k}(S_n) \}$ is infinite. Since $F_{k+1} \subseteq F_k$, the intersection of any finitely many sets from $\{F_k\}_{k\in\IN}$ is infinite, so the family $\{F_k\}_{k\in\IN}$ together with the Fr\'echet filter on $\IN$ has the finite intersection property. By the ultrafilter theorem (see, e.g., \cite[Theorem 1.1.6 and Corollary 1.1.18]{Goldbring2022}), this family extends to a non-principal ultrafilter $\omega \in \partial_\beta \IN$. For any $\delta > 0$ and any integer $N \ge 1$, choosing $k \ge \max\{N, 1/\delta\}$ gives $F_k \in \omega$ and $r_{1/k}(S_n) \ge r_\delta(S_n)$, so that
$$\rho_{\varepsilon,\omega}(T) > k\, \rho_{1/k,\omega}(S) \ge N\, \rho_{\delta,\omega}(S) \quad \text{in } \IN^\omega.$$
Thus $T \not\preceq_\omega S$, contradicting (3).
\end{proof}

Applying the faithful $*$-homomorphism $E: C^*_u(X)/\K \to \Q$ immediately yields the fundamental ultrarank obstruction for ideals in $C^*_u(X)$:

\begin{Cor}\label{cor: ultrarank obstruction}
Let $S, T \in C^*_u(X)$. If $T \in \langle S \rangle$, then $T \preceq S$. Consequently, if $T \not\preceq_\omega S$ for some $\omega \in \partial_\beta \IN$, then $T \notin \langle S \rangle$.

In particular, if $P = \bigoplus_{n\in\IN} P_n \in C^*_u(X)$ is a non-zero projection with $\sup_{n\in\IN} \rank(P_n) \le r < \infty$, then every $T \in \langle P \rangle$ satisfies $T \preceq P$, i.e., $\rho_{\varepsilon,\omega}(T)$ is finite in $\IN^\omega$ for all $\omega \in \partial_\beta \IN$.
\end{Cor}

\begin{proof}
If $T \in \langle S \rangle$, then $[E(T)] \in \langle [E(S)] \rangle_{\Q}$, so $T \preceq S$ follows directly from Theorem \ref{thm: ambient ideal ultrarank}.

When $P = \bigoplus_{n\in\IN} P_n$ is a non-zero projection with $\rank(P_n) \le r$, we have $r_\delta(P_n) \le r$ for all $\delta > 0$, and $r_\delta(P_n) \ge 1$ for $\delta \in (0,1)$ whenever $P_n \neq 0$. Thus $T \preceq P$ is equivalent to $\sup_{n\in\IN} r_\varepsilon(T_n) < \infty$ for every $\varepsilon > 0$ by Theorem \ref{thm: ambient ideal ultrarank} (2), and is equivalent to that $\rho_{\varepsilon,\omega}(T)$ is finite for all $\omega \in \partial_\beta \IN$ by Theorem \ref{thm: ambient ideal ultrarank} (3).
\end{proof}

\subsection{A counterexample via product expanders}

Before constructing the counterexample to \cite[Question 8.1]{WZ2025}, we record a diagonal cutoff lemma for geometric ideals, which follows directly from \cite{CW2004}.

\begin{Lem}\label{lem: diagonal cutoff}
Let $Y$ be a strongly discrete metric space with bounded geometry, and let $I \triangleleft C^*_u(Y)$ be a geometric ideal. For every $T \in I$ and every $\varepsilon > 0$, there exists a subset $B \subseteq Y$ such that $\chi_B \in I$ and
$$\|(1 - \chi_B)T\| < \varepsilon, \quad \|T(1 - \chi_B)\| < \varepsilon, \quad \text{and} \quad \|T - \chi_B T \chi_B\| < \varepsilon.$$
\end{Lem}

\begin{proof}
Since $I$ is geometric, we can choose a finite propagation operator $S \in I \cap \IC_u[Y]$ such that $\|T - S\| < \varepsilon/4$. For any $t > 0$, let $S_t$ be the $t$-truncation of $S$ defined by $S_t(x,y) = S(x,y)$ if $|S(x,y)| \ge t$ and $0$ otherwise. By \cite[Corollary 3.6]{CW2004}, $\lim_{t\to 0} \|S - S_t\| = 0$, so we may fix $t > 0$ such that $\|S - S_t\| < \varepsilon/4$. Let $B_r = r(\supp_t(S))$ and $B_c = r(\supp_t(S^*))$ be the row and column $t$-support sets of $S$. By \cite[Lemma 3.4 and Theorem 3.5]{CW2004}, both $\chi_{B_r}$ and $\chi_{B_c}$ belong to $\langle S \rangle_{C^*_u(Y)} \subseteq I$, and hence $B := B_r \cup B_c$ satisfies $\chi_B = \chi_{B_r} + \chi_{B_c} - \chi_{B_r}\chi_{B_c} \in I$. Since $S_t$ is supported on $B \times B$, we have $(1 - \chi_B)S_t = S_t(1 - \chi_B) = 0$. Consequently,
$$\|(1 - \chi_B)T\| \le \|T - S\| + \|S - S_t\| < \frac{\varepsilon}{2}, \quad \|T(1 - \chi_B)\| \le \|T - S\| + \|S - S_t\| < \frac{\varepsilon}{2},$$
and $\|T - \chi_B T \chi_B\| \le \|(1 - \chi_B)T\| + \|\chi_B T(1 - \chi_B)\| < \varepsilon$.
\end{proof}

Let $(E_n)_{n\in\IN}$ be a sequence of finite, connected, $d$-regular expander graphs with $\# E_n \to \infty$ as $n \to \infty$. Let $\Delta_n \in \B(\ell^2(E_n))$ be the Laplacian of $E_n$. By assumption, there exists a uniform spectral gap $\gamma > 0$ such that for all $n \in \IN$,
$$\sigma(\Delta_n) \subseteq \{0\} \cup [\gamma, 2d].$$
For each $n \in \IN$, set
$$X_n = E_n \times E_n, \quad \# X_n = (\# E_n)^2,$$
equipped with the $\ell^1$-product graph metric $d_{X_n}((a,b),(a',b')) = d_{E_n}(a,a') + d_{E_n}(b,b')$. We define $X = \bigsqcup_{n\in\IN} X_n$ to be the coarse disjoint union. Then $X$ is a strongly discrete metric space with bounded geometry.

Let $\xi_n = \frac{1}{\sqrt{\# E_n}} 1_{E_n} \in \ell^2(E_n)$ be the constant unit vector on $E_n$, and let $p_n = \xi_n \ox \xi_n \in \B(\ell^2(E_n))$ be the corresponding rank-one averaging projection. On $\ell^2(X_n) \cong \ell^2(E_n) \ox \ell^2(E_n)$, define
\begin{equation}\label{eq: P and Q}
P_n = p_n \ox p_n, \quad P = \bigoplus_{n\in\IN} P_n, \qquad P_{1,n} = p_n \ox 1, \quad P_1 = \bigoplus_{n\in\IN} P_{1,n}.
\end{equation}
Denote by $J_P = \langle P \rangle$ the closed two-sided ideal of $C^*_u(X)$ generated by $P$.

\begin{Thm}\label{thm: counterexample WZ}
For $X = \bigsqcup_{n\in\IN} (E_n \times E_n)$ constructed above, one has that
$$\K \subsetneq J_P \subsetneq \G \quad \text{and} \quad U(J_P) = X.$$
Moreover, there do not exist any geometric ideal $I \triangleleft C^*_u(X)$ such that
$$J_P = I \cap \G.$$
\end{Thm}

\begin{proof}
First, consider the bounded block-diagonal operators $D_1 = \bigoplus_{n\in\IN} (\Delta_n \ox 1)$ and $D_2 = \bigoplus_{n\in\IN} (1 \ox \Delta_n)$ on $\ell^2(X)$. Both $D_1$ and $D_2$ have propagation $1$ and norm at most $2d$, so $D_1, D_2 \in C^*_u(X)$. Choose a continuous function $f \in C([0, 2d])$ such that $f(0) = 1$ and $f|_{[\gamma, 2d]} \equiv 0$. Since $\ker(\Delta_n) = \IC \xi_n$, continuous functional calculus in $C^*_u(X)$ yields
$$P_1 = f(D_1) \in C^*_u(X) \quad \text{and} \quad P = f(D_1)f(D_2) \in C^*_u(X).$$
Within the $n$-th block $X_n$, the matrix entries of $P_n$ are all equal to $(\# E_n)^{-2}$, and the matrix entries of $P_{1,n}$ are equal to $(\# E_n)^{-1}$. Since $\# E_n \to \infty$, both $P$ and $P_1$ are non-compact ghost projections, so $P, P_1 \in \G \setminus \K$. Since $P \neq 0$ and $\K$ is the minimal non-zero ideal of $C^*_u(X)$, we have $\K \subsetneq J_P \subseteq \G$.

To see that $J_P \subsetneq \G$, notice that $\rank(P_n) = 1$ for every $n \in \IN$. By Corollary \ref{cor: ultrarank obstruction}, every operator $T \in J_P$ has finite $\varepsilon$-ultrarank $\rho_{\varepsilon,\omega}(T) \in \IN$ for all $\varepsilon > 0$ and $\omega \in \partial_\beta \IN$. On the other hand, $P_{1,n} = p_n \ox 1$ is a projection of rank $\# E_n$, so for every $\varepsilon \in (0,1)$,
$$\rho_{\varepsilon,\omega}(P_1) = [\# E_n]_\omega \in \IN^\omega \setminus \IN,$$
which is infinite. By Corollary \ref{cor: ultrarank obstruction}, $P_1 \notin J_P$, and hence $\K \subsetneq J_P \subsetneq \G$. Since $\K \cap \ell^\infty(X) = \G \cap \ell^\infty(X) = C_0(X)$, the inclusions $\K \subseteq J_P \subseteq \G$ imply $J_P \cap \ell^\infty(X) = C_0(X)$, which means $U(J_P) = X$ and $\wt I(U(J_P)) = \wt I(X) = \G$.

We now prove that $J_P \neq I \cap \G$ for any geometric ideal $I \triangleleft C^*_u(X)$. Assume for a contradiction that $J_P = I \cap \G$ for some geometric ideal $I$. Since $P \in J_P$, we have $P \in I$. Applying Lemma \ref{lem: diagonal cutoff} to $P \in I$ with $\varepsilon = 1/2$, we obtain a subset $B \subseteq X$ with $\chi_B \in I$ such that
$$\quad \|(1 - \chi_B)P\| < \frac{1}{2}.$$
Write $B_n = B \cap X_n$. Since $P_n = \xi_{X_n} \ox \xi_{X_n}$ where $\xi_{X_n} = \xi_n \ox \xi_n = (\# E_n)^{-1} 1_{X_n}$ is the constant unit vector on $X_n$, we have
$$\|(1 - \chi_{B_n})P_n\|^2 = \|(1 - \chi_{B_n})\xi_{X_n}\|^2 = 1 - \frac{\# B_n}{\# X_n} =1 - \frac{\# B_n}{(\# E_n)^2} < \frac{1}{4}.$$
Consequently,
\begin{equation}\label{eq: Bn density}
\# B_n > \frac{3}{4} (\# E_n)^2 \quad \text{for every } n \in \IN.
\end{equation}
Consider the restricted operator
$$R_B = \chi_B P_1 \chi_B.$$
Since $\chi_B \in I$, $P_1 \in \G$, we have $R_B \in I \cap \G = J_P$. Let us compute the $\varepsilon$-ultrarank of $R_B$. On each block $X_n$, the operator $(R_B)_n = \chi_{B_n}(p_n \ox 1)\chi_{B_n}$ decomposes as the orthogonal direct sum over $b \in E_n$ of the rank-one positive operators $\chi_{B_{n,b}} p_n \chi_{B_{n,b}}$, where $B_{n,b} = \{a \in E_n \mid (a,b) \in B_n\}$. The non-zero eigenvalue on the fiber of $b$ is
$$\lambda_{n,b} = \frac{\# B_{n,b}}{\# E_n} \in [0, 1].$$
By \eqref{eq: Bn density},
$$\sum_{b \in E_n} \lambda_{n,b} = \frac{\# B_n}{\# E_n} > \frac{3}{4} \# E_n.$$
If no more than $(\# E_n) / 2$ of the numbers $\{\lambda_{n,b}\}_{b \in E_n}$ were strictly greater than $1/2$, their sum is at most $\frac{\# E_n}{2} \cdot 1 + \frac{\# E_n}{2} \cdot \frac{1}{2} = \frac{3}{4} \# E_n$,  which leads to a contradiction. Therefore, at least $(\# E_n) / 2$ eigenvalues of $(R_B)_n$ are strictly greater than $1/2$, which means
$$r_{1/2}((R_B)_n) \ge \frac{\# E_n}{2} \quad \text{for all } n \in \IN, \quad \text{and} \quad \rho_{1/2,\omega}(R_B) \ge \Big[\frac{\# E_n}{2}\Big]_\omega.$$
Since $[(\# E_n)/2]_\omega$ is infinite in $\IN^\omega$, Corollary \ref{cor: ultrarank obstruction} implies that $R_B \notin J_P$. This is a contradiction.
\end{proof}

\section{Lower ultrarank does not imply ideal inclusion}\label{sec: converse failure}

In Theorem \ref{thm: ambient ideal ultrarank} and Corollary \ref{cor: ultrarank obstruction}, we saw that in the block-diagonal Calkin algebra $\Q$, the ultrarank order $T \preceq S$ is both necessary and sufficient for $[T] \in \langle [S] \rangle_{\Q}$, whereas in $C^*_u(X)$ it serves as a necessary obstruction to $T \in \langle S \rangle$. It is natural to ask whether the converse holds in $C^*_u(X)$, i.e., if $T \preceq S$, must $\langle S \rangle$ contain $T$?

In this section, we prove that the answer is negative. From a coarse geometric viewpoint, this phenomenon is actually very natural. Even when $S_n$ has much larger $\varepsilon$-rank than $T_n$ on $X_n$, the geometric regions where the singular vectors of $T_n$ and $S_n$ are concentrated may drift infinitely far apart as $n \to \infty$, so that no finite propagation operator in $C^*_u(X)$ can transport the rank of $S_n$ to approximate $T_n$.

Retain the notation of the product-expander space $X = \bigsqcup_{n\in\IN} X_n$ with $X_n = E_n \times E_n$ and $\# E_n \to \infty$ from Section \ref{sec: counterexample}. Consider the diagonal subset of $X_n$:
$$\Delta_n = \{ (a,a) \in X_n \mid a \in E_n \},$$
and let $\zeta_n = (\# E_n)^{-1/2} \sum_{a \in E_n} \delta_{(a,a)} \in \ell^2(X_n)$ be the normalized constant vector supported on $\Delta_n$. Define the block-diagonal rank-one projection
\begin{equation}\label{eq: R_Delta}
P_{\Delta,n} = \zeta_n \ox \zeta_n \in \B(\ell^2(X_n)), \quad P_\Delta = \bigoplus_{n\in\IN} P_{\Delta,n} \in \B(\ell^2(X)).
\end{equation}
Since $d_{X_n}((a,a),(b,b)) = 2 d_{E_n}(a,b)$ for all $a, b \in E_n$, the diagonal map $\iota: E \to X$ given by $a \mapsto (a,a)$ is an injective coarse embedding onto $\Delta = \bigsqcup_{n\in\IN} \Delta_n \subseteq X$. Consequently, the isometric embedding $V = \bigoplus_{n\in\IN} V_n: \ell^2(E) \to \ell^2(X)$ defined by $V_n \delta_a = \delta_{(a,a)}$ maps $\IC_u[E]$ into $\IC_u[X]$ by conjugation, and thus induces an isometric $*$-homomorphism $\textup{Ad}_V: C^*_u(E) \to C^*_u(X)$, $T \mapsto V T V^*$. Since $(E_n)_{n\in\IN}$ is an expander sequence, the averaging projection $p = \bigoplus_{n\in\IN} p_n$ is a non-compact ghost projection in $C^*_u(E)$. Noting that $V_n \xi_n = \zeta_n$, we immediately obtain that $P_\Delta = V p V^*$ lies in $C^*_u(X)$. 

For every $\varepsilon \in (0,1)$ and every non-principal ultrafilter $\omega \in \partial_\beta \IN$, it is direct that
$$\rho_{\varepsilon,\omega}(P_\Delta) = 1 \quad \text{and} \quad \rho_{\varepsilon,\omega}(P_1) = [\# E_n]_\omega,$$
so that $P_\Delta \simeq P \prec_\omega P_1$. 

\begin{Thm}\label{thm: low ultrarank no absorption}
For any $S\in \IC_u[X]$, one has that
$$P_1 S P_\Delta \in \K.$$
As a result, one has that
\begin{equation}\label{eq: intersection is K}
\langle P_\Delta \rangle \cap \langle P_1 \rangle = \K,
\end{equation}
and consequently $P_\Delta \notin \langle P_1 \rangle$.
\end{Thm}

\begin{proof}
Let $S \in \IC_u[X]$ be any finite propagation operator with $\prop(S) \le R$ and $\|S\| \le C$. Recall that $N_R = \sup_{n\in\IN, a\in E_n} \# B_{E_n}(a, R) < \infty$. On the $n$-th block $X_n$, the range of $P_{1,n} = p_n \ox 1$ has an orthonormal basis $\{\eta_b\}_{b \in E_n}$ given by
$$\eta_b = (\# E_n)^{-1/2} \sum_{a \in E_n} \delta_{(a,b)}, \quad b \in E_n.$$
For each $b \in E_n$, we compute the inner product
$$\langle \eta_b, S_n \zeta_n \rangle = \frac{1}{\# E_n} \sum_{a, c \in E_n} S_n\big((a,b), (c,c)\big).$$
Since $\prop(S_n) \le R$, a term $S_n((a,b),(c,c))$ in the sum can be non-zero only if
$$d_{X_n}\big((a,b),(c,c)\big) = d_{E_n}(a,c) + d_{E_n}(b,c) \le R.$$
For a fixed $b \in E_n$, there are at most $N_R$ choices for $c \in B_{E_n}(b, R)$, and for each such $c$, there are at most $N_R$ choices for $a \in B_{E_n}(c, R)$. Thus at most $N_R^2$ terms in the sum are non-zero, each of them has absolut value at most $\|S\| \le C$. It follows that for every $b \in E_n$,
$$|\langle \eta_b, S_n \zeta_n \rangle| \le \frac{N_R^2 C}{\# E_n} .$$
Summing the squares over the orthonormal basis $\{\eta_b\}_{b \in E_n}$ of $\textup{ran}(P_{1,n})$, we obtain
$$\|P_{1,n} S_n P_{\Delta,n}\| = \|P_{1,n} S_n \zeta_n\| = \Big( \sum_{b \in E_n} |\langle \eta_b, S_n \zeta_n \rangle|^2 \Big)^{1/2} \le \Big( \# E_n \cdot \frac{N_R^4 C^2}{(\# E_n)^2} \Big)^{1/2} = \frac{N_R^2 C}{\sqrt{\# E_n}} \to 0$$
as $n \to \infty$. This implies $P_1 E(S) P_\Delta \in \K$. Since $S - E(S) \in \K$, we have $P_1 S P_\Delta \in \K$. Approximating an arbitrary $T \in C^*_u(X)$ in norm by finite propagation operators, we have that
\begin{equation}\label{eq: Q A R_Delta in K}
P_1 C^*_u(X) P_\Delta \subseteq \K.
\end{equation}

From \eqref{eq: Q A R_Delta in K}, we have that $\langle P_1 \rangle \cdot \langle P_\Delta \rangle = \langle P_1 \rangle \cap \langle P_\Delta \rangle \subseteq \K$, which proves \eqref{eq: intersection is K}. Consequently $P_\Delta \notin \langle P_1 \rangle$.
\end{proof}

\begin{Rem}
Theorem \ref{thm: low ultrarank no absorption} shows that in $C^*_u(X)$, the finite propagation constraint creates geometric orthogonality at infinity between the fibers $E_n \times \{b\}$ and the diagonal $\Delta_n$. More precisely, any ball of radius $R$ in $X_n$ intersects both $E_n \times \{b\}$ and $\Delta_n$ in at most $N_R^2$ points. A similar phenomenon occurs also for a single fiber. If we pick a point $b_n \in E_n$ for each $n$ and set $B = \bigsqcup_{n\in\IN} (E_n \times \{b_n\})$, then the fiber ghost projection $P_{B} := \chi_B P_1 \chi_B = \bigoplus_{n\in\IN} (p_n \ox \chi_{\{b_n\}}) \in \G$ satisfies $P_{B} \simeq P$, yet $\chi_{B_n} P_n \chi_{B_n} = (\# E_n)^{-1} (p_n \ox \chi_{\{b_n\}})$ has norm $(\# E_n)^{-1} \to 0$, which shows $P_{B} \notin \langle P \rangle$.

Together with Section \ref{sec: counterexample}, these phenomena suggest that while the global ultrarank order $\preceq$ already distinguishes different ideals between $I(U)$ and $\wt I(U)$, a finer description of the intermediate ideal lattice $[I(U), \wt I(U)]$ should take into account \emph{local} ultrarank growth conditions, namely the growth rates of $\varepsilon$-ultraranks of operators compressed to coarse geometric subspaces of $X$.
\end{Rem}

\section{Ideal structure when $P$ generates the ghost ideal}\label{sec: positive results}

In this section, we prove that \cite[Question 8.1]{WZ2025} has an affirmative answer when the averaging (Kazhdan) projection $P$ generates the entire ghost ideal $\G$, and we establish a characterization of $\G = \langle P \rangle$ in terms of the $\varepsilon$-ultrarank and Yu's Property A on small-density subspaces.

Throughout this section, let $X = \bigsqcup_{n\in\IN} X_n$ be a coarse disjoint union of finite metric spaces with bounded geometry such that $\# X_n \to \infty$. Let
$$\xi_n = (\# X_n)^{-1/2} 1_{X_n} \in \H_n = \ell^2(X_n), \quad p_n = \xi_n \ox \xi_n \in \B(\H_n), \quad P = \bigoplus_{n\in\IN} p_n \in \B(\ell^2(X)).$$
By a theorem of A.~Khukhro, K.~Li, F.~Vigolo, and J.~Zhang \cite[Theorem 6.1]{KLVZ2021}, $P \in C^*_u(X)$ if and only if $(X_n)_{n\in\IN}$ is an \emph{asymptotic expander}. From now on, we assume that $P \in C^*_u(X)$ and denote $J_P = \langle P \rangle$.

\subsection{Characterization of $\langle P \rangle$ and zero-density Property A}

To characterize which operators belong to $J_P$, we first reduce the generators of $J_P$ to the form with multipliers in $\ell^\infty(X) \subseteq C^*_u(X)$.

\begin{Lem}\label{lem: normal form JP}
Assume that $P \in C^*_u(X)$. Then
$$J_P = \overline{\K + \textup{span}\{ f P g \mid f, g \in \ell^\infty(X) \}}^{\|\cdot\|}.$$
\end{Lem}

\begin{proof}
It is direct to see that the right-hand side is contained in $J_P$. Conversely, by \cite[Definition 2.3 and Lemma 2.4(2)]{CW2004} is spanned by operators of the form $f v$ (or $v g$), where $f, g \in \ell^\infty(X)$ and $v \in \IC_u[X]$ is a \emph{partial translation} induced by a bijection $t: D \to R$ between subsets $D, R \subseteq X$ with finite propagation. For all sufficiently large $n \in \IN$, $t$ restricts to a bijection from $D \cap X_n$ onto $R \cap X_n$. On each such block $X_n$, we have $v_n 1_{X_n} = \chi_{R \cap X_n} 1_{X_n}$ and hence $v_n p_n = \chi_{R \cap X_n} p_n$. It follows that $v P - \chi_R P \in \K$, and similarly $P v - P \chi_D \in \K$ by taking adjoints. Consequently, for any $A, B \in \IC_u[X]$, the operator $A P B$ belongs to $\K + \textup{span}\{ f P g \mid f, g \in \ell^\infty(X) \}$.
\end{proof}

For any $T \in C^*_u(X)$ and $\delta > 0$, define \emph{norm-density function}
\begin{equation}\label{eq: a_T delta}
a_T(\delta) := \sup_{\substack{n \ge 1,\ B \subseteq X_n \\ \# B \le \delta \# X_n}} \max\big\{ \|\chi_B T_n\|,\, \|T_n \chi_B\| \big\}.
\end{equation}
We say that $T$ is \emph{uniformly continuous in density}, if
$$\lim_{\delta \to 0} a_T(\delta) = 0.$$

\begin{Pro}\label{thm: characterization of JP}
Assume that $P \in C^*_u(X)$. An operator $T \in C^*_u(X)$ satisfies $T \in J_P$ if and only if it satisfies the following two conditions:
\begin{itemize}
\item[(1)] $T$ is uniformly continuous in density, and
\item[(2)] $T \preceq P$. 
\end{itemize}
\end{Pro}

We remark that the condition $T \preceq P$ is equivalent to for every $\varepsilon > 0$ and every $\omega \in \partial_\beta \IN$, the $\varepsilon$-ultrarank $\rho_{\varepsilon,\omega}(T)$ is finite in $\IN^\omega$, or equivalently, $\sup_{n\in\IN} r_\varepsilon(T_n) < \infty$, by Corollary \ref{cor: ultrarank obstruction}.

\begin{proof}
$( \Rightarrow )$: If $T \in J_P$, then $T \preceq P$ (and hence $\sup_{n\in\IN} r_\varepsilon(T_n) < \infty$ for every $\varepsilon > 0$) follows directly from Corollary \ref{cor: ultrarank obstruction}. To see that $T$ is uniformly continuous in density, note that for any $f, g \in \ell^\infty(X)$ and any subset $B \subseteq X_n$,
$$\|\chi_B f p_n g\| \le \|f\|\cdot \|g\| \cdot \|\chi_B \xi_n\| = \|f\| \cdot \|g\| \cdot \sqrt{\frac{\# B}{\# X_n}},$$
and the same estimate holds for $\|f p_n g \chi_B\|$. Thus any finite sum $U = \sum_{i=1}^k f_i P g_i$ satisfies $a_U(\delta) \le C_U \sqrt{\delta} \to 0$ as $\delta \to 0$. Moreover, any finitely supported block-diagonal operator $C = \bigoplus_{n=1}^N C_n \in \K$ is uniformly continuous in density, since on finitely many blocks, the condition $\# B \le \delta \# X_n$ forces $B = \emptyset$ when $\delta$ is sufficiently small. Since $a_{T - S}(\delta) \le \|T - S\|$, the set of operators that are uniformly continuous in density is norm-closed, so Lemma \ref{lem: normal form JP} implies that every $T \in J_P$ is uniformly continuous in density.

$( \Leftarrow )$: Note that on each block $X_n$, a rank-one operator of the form $u_n \ox v_n$ can be written as $f_n p_n g_n = (f_n \xi_n) \ox (\overline{g_n} \xi_n)$ by setting $f_n(x) = \sqrt{\# X_n}\, u_n(x)$ and $g_n(x) = \sqrt{\# X_n}\, \overline{v_n(x)}$. Moreover, the resulting functions $f = (f_n)$ and $g = (g_n)$ belong to $\ell^\infty(X)$ if and only if there exists a constant $D > 0$ independent of $n$ such that $|u_n(x)|, |v_n(x)| \le D / \sqrt{\# X_n}$ for all $x \in X_n$. Thus, to show that $T \in J_P$, our strategy is to first truncate the SVD of $T_n$ to a uniformly bounded number of rank-one terms using (2), and then truncate the coordinates of the corresponding singular vectors at height $D / \sqrt{\# X_n}$ using (1).

Fix $\eta > 0$. For each $n \in \IN$, write the SVD of $T_n$ as
$$T_n = \sum_{m=1}^{\# X_n} \tau_n^{(m)} \big(u_n^{(m)} \ox v_n^{(m)}\big),$$
where $\tau_n^{(1)} \ge \cdots \ge \tau_n^{(\# X_n)} \ge 0$ are the singular values of $T_n$, and $\{u_n^{(m)}\}_{m=1}^{\# X_n}$, $\{v_n^{(m)}\}_{m=1}^{\# X_n}$ are orthonormal bases of $\H_n$. Set $q_n = r_\eta(T_n)$, and let
$$T_n^\eta = \sum_{m=1}^{q_n} \tau_n^{(m)} \big(u_n^{(m)} \ox v_n^{(m)}\big)$$
be the SVD truncation of $T_n$ at $\eta$, so that $\|T_n - T_n^\eta\| \le \eta$. By (2), we have $K_\eta := \sup_{n\in\IN} q_n = \sup_{n\in\IN} r_\eta(T_n) < \infty$. For convenience, set $\tau_n^{(m)} = 0$ and $u_n^{(m)} = v_n^{(m)} = 0$ whenever $m > q_n$.

Given $D > 0$, to bound the coordinates of $u_n^{(m)}$ and $v_n^{(m)}$ by $D / \sqrt{\# X_n}$, we consider
$$B_n^{(m), D} = \Big\{ x \in X_n \ \Big|\ |u_n^{(m)}(x)| > \frac{D}{\sqrt{\# X_n}} \Big\}, \qquad \wt B_n^{(m), D} = \Big\{ x \in X_n \ \Big|\ |v_n^{(m)}(x)| > \frac{D}{\sqrt{\# X_n}} \Big\},$$
and define the coordinate-truncated vectors in $\H_n$ by
$$u_n^{(m), D} = \big(1 - \chi_{B_n^{(m), D}}\big) u_n^{(m)}, \qquad v_n^{(m), D} = \big(1 - \chi_{\wt B_n^{(m), D}}\big) v_n^{(m)}.$$
By construction, for each $1 \le m \le K_\eta$, the functions $f_m, g_m$ defined on each block $X_n$ by
$$f_m(x) = \tau_n^{(m)} \sqrt{\# X_n}\, u_n^{(m), D}(x) \quad \text{and} \quad g_m(x) = \sqrt{\# X_n}\, \overline{v_n^{(m), D}(x)} \qquad (x \in X_n)$$
satisfy $\|f_m\|_\infty \le \|T\| D$ and $\|g_m\|_\infty \le D$, so that $W^D := \bigoplus_{n\in\IN} W_n^D = \sum_{m=1}^{K_\eta} f_m P g_m \in J_P$, where $W_n^D = \sum_{m=1}^{K_\eta} \tau_n^{(m)} \big(u_n^{(m), D} \ox v_n^{(m), D}\big)$.

It remains to estimate the truncation error $\|T_n^\eta - W_n^D\|$. For any subset $B \subseteq X_n$ and $1 \le m \le q_n$, since $T_n v_n^{(m)} = \tau_n^{(m)} u_n^{(m)}$ and $T_n^* u_n^{(m)} = \tau_n^{(m)} v_n^{(m)}$ with $\tau_n^{(m)} > \eta$, we have
\begin{equation}\label{eq: singular vector tail}
\|\chi_B u_n^{(m)}\| = \frac{\|\chi_B T_n v_n^{(m)}\|}{\tau_n^{(m)}} \le \frac{\|\chi_B T_n\|}{\eta}, \qquad \|\chi_B v_n^{(m)}\| = \frac{\|\chi_B T_n^* u_n^{(m)}\|}{\tau_n^{(m)}} \le \frac{\|T_n \chi_B\|}{\eta}.
\end{equation}
On the other hand, since $\|u_n^{(m)}\| = 1$ and $\|v_n^{(m)}\| = 1$, Chebyshev's inequality implies that $\# B_n^{(m), D} \le D^{-2} \# X_n$ and $\# \wt B_n^{(m), D} \le D^{-2} \# X_n$. Combining \eqref{eq: singular vector tail} with \eqref{eq: a_T delta} yields
$$\|u_n^{(m)} - u_n^{(m), D}\| = \|\chi_{B_n^{(m), D}} u_n^{(m)}\| \le \frac{a_T(D^{-2})}{\eta}, \qquad \|v_n^{(m)} - v_n^{(m), D}\| = \|\chi_{\wt B_n^{(m), D}} v_n^{(m)}\| \le \frac{a_T(D^{-2})}{\eta}.$$
Consequently, Since $T$ is uniformly continuous in density, we have
\begin{equation}\label{eq: Vn Wn approx}
\sup_{n\in\IN} \|T_n^\eta - W_n^D\| \le \frac{2 \|T\| K_\eta}{\eta}\, a_T(D^{-2}) \to 0 \quad \text{as } D \to \infty.
\end{equation}
Let $D \to \infty$ and then $\eta \to 0$, we conclude that $T \in J_P$.
\end{proof}

To see the coarse geometric meaning of uniform continuity in density, observe that if $T \in \G$ is uniformly continuous in density, then for any sequence of subsets $Y_n \subseteq X_n$ with $\lim_{n\to\infty} \frac{\# Y_n}{\# X_n} = 0$, the norms $\|\chi_{Y_n} T_n\|$ and $\|T_n \chi_{Y_n}\|$ must converge to $0$ as $n \to \infty$, which implies that the compression of $T$ to the zero-density subspace $Y = \bigsqcup_{n\in\IN} Y_n$ is compact. Thus, requiring all ghost operators on $X$ to be uniformly continuous in density means precisely that no zero-density subspace $Y \subseteq X$ can support a non-compact ghost operator (in contrast to the diagonal subspace $\Delta$ in Theorem \ref{thm: low ultrarank no absorption}). By the Roe-Willett theorem \cite[Theorem 1.3]{RoeWil2014}, the absence of non-compact ghost operators on $Y$ is equivalent to $Y$ having Yu's Property A. This naturally leads to the following definition.

\begin{Def}\label{def: A0}
Let $X = \bigsqcup_{n\in\IN} X_n$. We say that $X$ has \emph{zero-density Property A} if for every sequence of subsets $Y_n \subseteq X_n$ satisfying $\lim_{n\to\infty} \frac{\# Y_n}{\# X_n} = 0$, the induced metric subspace $Y = \bigsqcup_{n\in\IN} Y_n \subseteq X$ has Yu's Property A.
\end{Def}

\begin{Pro}\label{thm: A0 iff AC}
Let $X = \bigsqcup_{n\in\IN} X_n$ with $\# X_n \to \infty$. The following are equivalent:
\begin{itemize}
\item[(1)] $X$ has zero-density Property A;
\item[(2)] every ghost operator $T \in \G$ is uniformly continuous in density.
\end{itemize}
\end{Pro}

\begin{proof}
$(1) \Rightarrow (2)$: Suppose $X$ has zero-density Property A. Let $Y_n \subseteq X_n$ be any sequence of subsets with $\frac{\# Y_n}{\# X_n} \to 0$, and set $Y = \bigsqcup_{n\in\IN} Y_n$. The uniform Roe algebra $C^*_u(Y)$ identifies with the corner $\chi_Y C^*_u(X) \chi_Y$. For any $T \in \G$, the operators $\chi_Y T T^* \chi_Y$ and $\chi_Y T^* T \chi_Y$ are ghost operators in $C^*_u(Y)$. Since $Y$ has Property A since $Y$ has zero-density, \cite[Theorem 1.3]{RoeWil2014} implies that $\chi_Y T T^* \chi_Y$ and $\chi_Y T^* T \chi_Y$ are compact. Thus, both $e_Y T$ and $T e_Y$ are compact, i.e.,
\begin{equation}\label{eq: Yn compression to 0}
\lim_{n\to\infty} \max\big\{ \|\chi_{Y_n} T_n\|,\, \|T_n \chi_{Y_n}\| \big\} = 0.
\end{equation}
However, if $T$ is not uniformly continuous in density, then there exists $\varepsilon_0 > 0$ such that for any $k\in\IN$, there exists a subset $B_{n_k} \subseteq X_{n_k}$ with $\# B_{n_k} \le \frac{1}{k} \# X_{n_k}$ such that 
$$\max\{\|\chi_{B_{n_k}} T_{n_k}\|, \|T_{n_k} \chi_{B_{n_k}}\|\} > \varepsilon_0.$$
Setting $Y_{n_k} = B_{n_k}$ and $Y_n = \emptyset$ for $n \notin \{n_k\}$, this leads to a contradiction with \eqref{eq: Yn compression to 0}.

$(2) \Rightarrow (1)$: Assume that every $T \in \G$ is uniformly continuous in density. Let $Y = \bigsqcup_{n\in\IN} Y_n$ be any subspace with $\frac{\# Y_n}{\# X_n} \to 0$, and let $T \in C^*_u(Y)\subseteq C^*_u(X)$ be a ghost operator on $Y$. Given any $\varepsilon > 0$, there exists $\delta > 0$ such that
$$\max\{\|\chi_B T_n\|, \|T_n \chi_B\|\} < \varepsilon$$
whenever $\# B \le \delta \# X_n$. For all sufficiently large $n$, we have $\# Y_n \le \delta \# X_n$. Thus, we get $\|T_n\| = \|\chi_{Y_n} T_n\| < \varepsilon$. Therefore, $\|T_n\| \to 0$ as $n \to \infty$, which implies $T \in \K(\ell^2(Y))$. By \cite[Theorem 1.3]{RoeWil2014}, $Y$ has Property A.
\end{proof}

Combining Proposition \ref{thm: characterization of JP} and Proposition \ref{thm: A0 iff AC}, we immediately obtain the characterization of when the averaging projection $P$ generates the entire ghost ideal $\G$:

\begin{Thm}\label{thm: G=JP characterization}
Let $X = \bigsqcup_{n\in\IN} X_n$ be a sequence of asymptotic expander graphs. Then $\G = \langle P \rangle$ if and only if both of the following conditions hold:
\begin{itemize}
\item[(1)] for every ghost operator $T \in \G$, one has $T \preceq P$;
\item[(2)] $X$ has zero-density Property A.\qed
\end{itemize}
\end{Thm}

\subsection{Classification of ideals when $\G = \langle P \rangle$}

We now show that when $\G = \langle P \rangle$, \cite[Question 8.1]{WZ2025} has an affirmative answer. To extend the classification from the interval $[\K, \G]$ to all closed two-sided ideals of $C^*_u(X)$, we first show that $\langle P \rangle = \G$ yields the decomposition $\wt I(U) = I(U) + \G$ for every non-empty invariant open subset $U \subseteq \beta X$.

\begin{Lem}\label{lem: ghostly decomposition}
Assume that $P \in C^*_u(X)$ and $\G = \langle P \rangle$. Then for every non-empty $G(X)$-invariant open subset $U \subseteq \beta X$, one has that
$$\wt I(U) = I(U) + \G.$$
\end{Lem}

\begin{proof}
Let $U \subseteq \beta X$ be a non-empty $G(X)$-invariant open subset. Since $X$ is the smallest invariant open set which is contained in any invariant open set $U$, we have that $\G \subseteq \wt I(U)$. Consequently, $I(U) + \G \subseteq \wt I(U)$.

Conversely, let $T \in \wt I(U)$. Without loss of generality, we assume that $T = \bigoplus_{n\in\IN} T_n$ is block-diagonal. For each $k \in \IN$, consider the $1/k$-support of $T$
$$B_k = \big\{ x \in X \ \big|\ \exists y \in X \text{ such that } |T(x,y)| \ge 1/k \text{ or } |T(y,x)| \ge 1/k \big\},$$
and define the remainder
$$R_k = (1 - \chi_{B_k}) T (1 - \chi_{B_k}) = \bigoplus_{n\in\IN} (R_k)_n.$$
Since $T \in \wt I(U)$, we have $\overline{B_k}^{\beta X} \subseteq U$ and hence $\chi_{B_k} \in I(U)$. It follows that $T - R_k = \chi_{B_k} T + (1 - \chi_{B_k}) T \chi_{B_k} \in I(U)$. Since $I(U) + \G$ is closed in $C^*_u(X)$, to prove that $T \in I(U) + \G$, we shall now prove that
$$\inf_{k\in\IN} d(R_k, \G) = 0.$$

Since $\G = \langle P \rangle$, for each $m \in \IN$, define
$$C_m = \left\{ \sum_{j=1}^m f_j P g_j \ \Big|\ f_j, g_j \in \ell^\infty(X),\ \|f_j\|, \|g_j\| \le m \right\},$$
and let $C_{m,n} = \{ \chi_{X_n} W \chi_{X_n} \mid W \in C_m \} \subseteq \B(\H_n)$ be its $n$-th block. Since modifying finitely many blocks only changes an operator by an element of $\K$, for any block-diagonal operator $S = \bigoplus_{n\in\IN} S_n \in C^*_u(X)$ we have
$$d(S, C_m + \K) = \limsup_{n\to\infty} d\big(S_n, C_{m,n}\big).$$
By Lemma \ref{lem: normal form JP}, $\G = \langle P \rangle$ is the norm closure of $\bigcup_{m\in\IN} (C_m + \K)$, which yields
\begin{equation}\label{eq: dist formula G}
d(S, \G) = \lim_{m\to\infty} \limsup_{n\to\infty} d\big(S_n, C_{m,n}\big).
\end{equation}

Apply \eqref{eq: dist formula G} to $R_k$ and take $m = k$, we obtain
$$d(R_k, \G) \le \limsup_{n\to\infty} d\big((R_k)_n, C_{k,n}\big) \quad \text{for all } k \in \IN.$$
Assume for a contradiction that there exists $\varepsilon_0 > 0$ such that $d(R_k, \G) \ge \varepsilon_0$ for all $k \in \IN$. Then
$$\limsup_{n\to\infty} d((R_k)_n, C_{k,n}) \geq d(R_k, \G) > \varepsilon_0$$
for every $k \in \IN$. We could then inductively choose a strictly increasing sequence $(n_k)_{k\in\IN}$ in $\IN$ such that
$$d\big((R_k)_{n_k}, C_{k,n_k}\big) > \varepsilon_0 \quad \text{for all } k \in \IN.$$
Let $f \in \ell^\infty(X) \subseteq C^*_u(X)$ be the diagonal projection defined by $f|_{X_{n_k}} = \chi_{X_{n_k} \setminus B_k}$ for each $k \in \IN$ and $f|_{X_n} = 0$ otherwise, and set $S = f T f \in C^*_u(X)$. Then $S$ satisfies that $S_{n_k} = (R_k)_{n_k}$ and $S_n = 0$ for $n \notin \{n_k\}_{k\in\IN}$. Since all entries of $S_{n_k} = (R_k)_{n_k}$ have absolute value less than $1/k$, $S$ is a ghost operator in $\G$. However, for any fixed $m \in \IN$, if $k>m$, we have that
$$d(S_{n_k}, C_{m,n_k}) \ge d(S_{n_k}, C_{k,n_k}) > \varepsilon_0,$$
since $C_{m,n_k} \subseteq C_{k,n_k}$. Thus, for any fixed $m \in \IN$, we have
$$d(S, C_m + \K) = \limsup_{n\to\infty} d\big(S_n, C_{m,n}\big) > \varepsilon_0.$$
By \eqref{eq: dist formula G}, this means that $d(S, \G) > \varepsilon_0$, which is a contradiction.
\end{proof}

\begin{Thm}\label{thm: full classification under G=JP}
Assume that $P \in C^*_u(X)$, and $\G = \langle P \rangle$. Then for every closed two-sided ideal $J \in [\K, \G]$ in $C^*_u(X)$, there exists a $G(X)$-invariant open subset $V \subseteq \beta X$ such that
$$ J = I(V) \cap \G. $$
\end{Thm}

\begin{proof}
Let $J \in [\K, \G]$. Since $\G = \langle P \rangle$ and $J \subseteq \G$, we have
$$J = \G J \G = \langle P J P \rangle.$$
For any $T \in J$, without loss of generality, we shall assume $T$ is block diagonal. It is a direct computation that on each block $X_n$, $p_n T_n p_n = \lambda_n p_n$ where $\lambda_n = \langle T_n \xi_n, \xi_n \rangle$. Define $f_T \in \ell^\infty(X)$ to be the function given by $f_T|_{X_n} = \lambda_n \chi_{X_n}$ for each $n \in \IN$, so that $P T P = f_T P \in J$. For any $\varepsilon > 0$, let
$$B_\varepsilon = \big\{x \in X \mid |f_T(x)| \ge \varepsilon\big\} = \bigsqcup_{n\text{ such that }|\lambda_n| \ge \varepsilon} X_n,$$
and define $h_\varepsilon \in \ell^\infty(X)$ by $h_\varepsilon(x) = f_T(x)^{-1}$ on $B_\varepsilon$ and $0$ elsewhere. Then $\chi_{B_\varepsilon} P = h_\varepsilon f_T P \in J$. Since $\|f_T P - f_T \chi_{B_\varepsilon} P\| \le \varepsilon$, set
$$\mathcal{F}_J = \Big\{ B \subseteq X \ \Big|\ B = \bigsqcup_{n \in S} X_n \text{ for some } S \subseteq \IN,\text{ and } \chi_B P \in J \Big\},$$
we have $f_T P \in \langle \{\chi_B P \mid B \in \mathcal{F}_J\} \rangle$ for every $T \in J$. Hence
$$J = \big\langle \{ \chi_B P \mid B \in \mathcal{F}_J \} \big\rangle.$$

Now, let $V=\bigcup_{B \in \mathcal{F}_J} \overline{B}^{\beta X}$. Then $V$ is a $G(X)$-invariant open subset of $\beta X$ and
$$I(V) = \langle \{ \chi_B \mid B \in \mathcal{F}_J \} \rangle$$
is a geometric ideal. From the generators, we directly see that $J \subseteq I(V) \cap \G$. Conversely, for each $B \in \mathcal{F}_J$, since $B$ is a union of blocks $X_n$ and $d(X_n, X \setminus X_n) \to \infty$, any finite propagation operator connects $B$ and $X \setminus B$ on at most finitely many blocks, which implies $[\chi_B, T] \in \K$ for all $T \in C^*_u(X)$. Thus, for any $T, S \in C^*_u(X)$,
$$\chi_B T P S = T (\chi_B P) S + [\chi_B, T] P S \in J + \K = J.$$
Since $\G = \langle P \rangle$, it follows that $\chi_B \G \subseteq J$ for all $B \in \mathcal{F}_J$, and hence $I(V) \G = I(V) \cap \G \subseteq J$. To sum up, we obtain $J = I(V) \cap \G$.
\end{proof}

\begin{Cor}\label{cor: general ideal classification}
Assume that $P \in C^*_u(X)$, and $\G = \langle P \rangle$. Then for every nonzero closed two-sided ideal $J \triangleleft C^*_u(X)$, there exists a $G(X)$-invariant open subset $V \supseteq U(J)$ in $\beta X$ such that
$$J = I(V) \cap \wt I(U(J)).$$
\end{Cor}

\begin{proof}
Let $J \triangleleft C^*_u(X)$ be a closed two-sided ideal, and set $U = U(J)$. By Lemma \ref{lem: ghostly decomposition}, we have $\wt I(U) = I(U) + \G$. Since $I(U) \subseteq J \subseteq \wt I(U) = I(U) + \G$, we can apply the distributive law of ideals to obtain
$$J = J \cap (I(U) + \G) = I(U) + (J \cap \G).$$
The intersection $J \cap \G$ belongs to $[\K, \G]$. By Theorem \ref{thm: full classification under G=JP}, there exists a $G(X)$-invariant open subset $V_0 \subseteq \beta X$ such that $J \cap \G = I(V_0) \cap \G$. Setting $V := U \cup V_0$, which is a $G(X)$-invariant open subset containing $U$, we have $I(V) = I(U) + I(V_0)$. Applying the distributive law once more yields
$$\begin{aligned}
J &= I(U) + (I(V_0) \cap \G) \\
&= (I(U) + I(V_0)) \cap (I(U) + \G) \\
&= I(V) \cap \wt I(U(J)).
\end{aligned}$$
This finishes the proof.
\end{proof}

\bibliographystyle{alpha}
\bibliography{ref}

@book {Roe2003,
    AUTHOR = {Roe, J.},
     TITLE = {Lectures on coarse geometry},
    SERIES = {University Lecture Series},
    VOLUME = {31},
 PUBLISHER = {American Mathematical Society, Providence, RI},
      YEAR = {2003},
     PAGES = {viii+175},
      ISBN = {0-8218-3332-4},
   MRCLASS = {53C24 (20F65 46L05 51K05 54E25)},
  MRNUMBER = {2007488},
MRREVIEWER = {Jean-Louis\ Tu},
       DOI = {10.1090/ulect/031},
       URL = {https://doi.org/10.1090/ulect/031},
}

@article {BBFKVW2022,
    AUTHOR = {Baudier, F. P. and Braga, B. M. and Farah, I. and
              Khukhro, A. and Vignati, A. and Willett, R.},
     TITLE = {Uniform {R}oe algebras of uniformly locally finite metric
              spaces are rigid},
   JOURNAL = {Invent. Math.},
  FJOURNAL = {Inventiones Mathematicae},
    VOLUME = {230},
      YEAR = {2022},
    NUMBER = {3},
     PAGES = {1071--1100},
      ISSN = {0020-9910,1432-1297},
   MRCLASS = {46L05 (51F30)},
  MRNUMBER = {4506772},
       DOI = {10.1007/s00222-022-01140-x},
       URL = {https://doi.org/10.1007/s00222-022-01140-x},
}

@article {ChenWang2005,
    AUTHOR = {Chen, X. and Wang, Q.},
     TITLE = {Ghost ideals in uniform {R}oe algebras of coarse spaces},
   JOURNAL = {Arch. Math. (Basel)},
  FJOURNAL = {Archiv der Mathematik},
    VOLUME = {84},
      YEAR = {2005},
    NUMBER = {6},
     PAGES = {519--526},
      ISSN = {0003-889X,1420-8938},
   MRCLASS = {46L80},
  MRNUMBER = {2148492},
MRREVIEWER = {Jean-Louis\ Tu},
       DOI = {10.1007/s00013-005-1189-1},
       URL = {https://doi.org/10.1007/s00013-005-1189-1},
}

@article {CW2004,
    AUTHOR = {Chen, X. and Wang, Q.},
     TITLE = {Ideal structure of uniform {R}oe algebras of coarse spaces},
   JOURNAL = {J. Funct. Anal.},
  FJOURNAL = {Journal of Functional Analysis},
    VOLUME = {216},
      YEAR = {2004},
    NUMBER = {1},
     PAGES = {191--211},
      ISSN = {0022-1236,1096-0783},
   MRCLASS = {46L80 (46L05)},
  MRNUMBER = {2091361},
MRREVIEWER = {Jean-Louis\ Tu},
       DOI = {10.1016/j.jfa.2003.11.015},
       URL = {https://doi.org/10.1016/j.jfa.2003.11.015},
}

@book {Goldbring2022,
    AUTHOR = {Goldbring, I.},
     TITLE = {Ultrafilters throughout mathematics},
    SERIES = {Graduate Studies in Mathematics},
    VOLUME = {220},
 PUBLISHER = {American Mathematical Society, Providence, RI},
      YEAR = {2022},
     PAGES = {xviii+399},
      ISBN = {978-1-4704-6900-9},
   MRCLASS = {03-02 (03C20 03C50 03E55 03H05 54D80)},
  MRNUMBER = {4454845},
MRREVIEWER = {Klaas\ Pieter\ Hart},
       DOI = {10.1090/gsm/220},
       URL = {https://doi.org/10.1090/gsm/220},
}

@book {HJ2013,
    AUTHOR = {Horn, R. A. and Johnson, C. R.},
     TITLE = {Matrix analysis},
   EDITION = {Second},
 PUBLISHER = {Cambridge University Press, Cambridge},
      YEAR = {2013},
     PAGES = {xviii+643},
      ISBN = {978-0-521-54823-6},
   MRCLASS = {15-01},
  MRNUMBER = {2978290},
MRREVIEWER = {Mohammad\ Sal\ Moslehian},
}

@article {KLVZ2021,
    AUTHOR = {Khukhro, A. and Li, K. and Vigolo, F. and Zhang, J.},
     TITLE = {On the structure of asymptotic expanders},
   JOURNAL = {Adv. Math.},
  FJOURNAL = {Advances in Mathematics},
    VOLUME = {393},
      YEAR = {2021},
     PAGES = {Paper No. 108073, 35},
      ISSN = {0001-8708,1090-2082},
   MRCLASS = {46H35 (05C99 19K56 20F65 51F30)},
  MRNUMBER = {4340225},
MRREVIEWER = {Hang\ Wang},
       DOI = {10.1016/j.aim.2021.108073},
       URL = {https://doi.org/10.1016/j.aim.2021.108073},
}

@book {RRS2004,
    AUTHOR = {Rabinovich, V. and Roch, S. and Silbermann, B.},
     TITLE = {Limit operators and their applications in operator theory},
    SERIES = {Operator Theory: Advances and Applications},
    VOLUME = {150},
 PUBLISHER = {Birkh\"auser Verlag, Basel},
      YEAR = {2004},
     PAGES = {xvi+392},
      ISBN = {3-7643-7081-5},
   MRCLASS = {47-02 (35-02 35S05 47A50 47A53 47B35 47G30)},
  MRNUMBER = {2075882},
MRREVIEWER = {Oleksiy\ Karlovych},
       DOI = {10.1007/978-3-0348-7911-8},
       URL = {https://doi.org/10.1007/978-3-0348-7911-8},
}

@article {RoeWil2014,
    AUTHOR = {Roe, J. and Willett, R.},
     TITLE = {Ghostbusting and property {A}},
   JOURNAL = {J. Funct. Anal.},
  FJOURNAL = {Journal of Functional Analysis},
    VOLUME = {266},
      YEAR = {2014},
    NUMBER = {3},
     PAGES = {1674--1684},
      ISSN = {0022-1236,1096-0783},
   MRCLASS = {46L70 (20F67 51F99)},
  MRNUMBER = {3146831},
MRREVIEWER = {J\'an\ \v Spakula},
       DOI = {10.1016/j.jfa.2013.07.004},
       URL = {https://doi.org/10.1016/j.jfa.2013.07.004},
}

@article {STY2002,
    AUTHOR = {Skandalis, G. and Tu, J. L. and Yu, G.},
     TITLE = {The coarse {B}aum-{C}onnes conjecture and groupoids},
   JOURNAL = {Topology},
  FJOURNAL = {Topology. An International Journal of Mathematics},
    VOLUME = {41},
      YEAR = {2002},
    NUMBER = {4},
     PAGES = {807--834},
      ISSN = {0040-9383},
   MRCLASS = {58J22 (19K56 46L80 46L85)},
  MRNUMBER = {1905840},
MRREVIEWER = {Paul\ D.\ Mitchener},
       DOI = {10.1016/S0040-9383(01)00004-0},
       URL = {https://doi.org/10.1016/S0040-9383(01)00004-0},
}

@article {SpaWil2017,
    AUTHOR = {{\v S}pakula, J. and Willett, R.},
     TITLE = {A metric approach to limit operators},
   JOURNAL = {Trans. Amer. Math. Soc.},
  FJOURNAL = {Transactions of the American Mathematical Society},
    VOLUME = {369},
      YEAR = {2017},
    NUMBER = {1},
     PAGES = {263--308},
      ISSN = {0002-9947,1088-6850},
   MRCLASS = {47A53 (30L05 46L85 47B36)},
  MRNUMBER = {3557774},
MRREVIEWER = {Ivan\ G.\ Todorov},
       DOI = {10.1090/tran/6660},
       URL = {https://doi.org/10.1090/tran/6660},
}

@book {WilYu2020,
    AUTHOR = {Willett, R. and Yu, G.},
     TITLE = {Higher index theory},
    SERIES = {Cambridge Studies in Advanced Mathematics},
    VOLUME = {189},
 PUBLISHER = {Cambridge University Press, Cambridge},
      YEAR = {2020},
     PAGES = {xi+581},
      ISBN = {978-1-108-49106-8},
   MRCLASS = {46-02 (19Kxx 46L80 47L60 58J22)},
  MRNUMBER = {4411373},
MRREVIEWER = {Man-Ho\ Ho},
       DOI = {10.1017/9781108867351},
       URL = {https://doi.org/10.1017/9781108867351},
}

@article {WZ2025,
    AUTHOR = {Wang, Q. and Zhang, J.},
     TITLE = {Ideal structure of uniform {R}oe algebras: beyond {P}roperty {A}},
   JOURNAL = {Proc. Roy. Soc. Edinburgh Sect. A},
  FJOURNAL = {Proceedings of the Royal Society of Edinburgh. Section A.
              Mathematics},
      YEAR = {2025},
     PAGES = {1--39},
       DOI = {10.1017/prm.2025.10095},
       URL = {https://doi.org/10.1017/prm.2025.10095},
}

@article {Yu2000,
    AUTHOR = {Yu, G.},
     TITLE = {The coarse {B}aum-{C}onnes conjecture for spaces which admit a
              uniform embedding into {H}ilbert space},
   JOURNAL = {Invent. Math.},
  FJOURNAL = {Inventiones Mathematicae},
    VOLUME = {139},
      YEAR = {2000},
    NUMBER = {1},
     PAGES = {201--240},
      ISSN = {0020-9910,1432-1297},
   MRCLASS = {19K56 (46L80 57R67 58J22)},
  MRNUMBER = {1728880},
MRREVIEWER = {Nigel\ Higson},
       DOI = {10.1007/s002229900032},
       URL = {https://doi.org/10.1007/s002229900032},
}

@article{GQW2024,
  author  = {Guo, L. and Qian, J. and Wang, Q.},
  title   = {A nonstandard analysis approach to limit operators and {Fredholmness} in {Roe}-like algebras},
  journal = {\textup{To appear in} Bull. Lond. Math. Soc.},
  year    = {2026},
  note    = {\href{https://arxiv.org/abs/2412.08130}{arXiv:2412.08130}}
}

\end{document}